\documentclass[12pt]{article}
\usepackage{amsfonts,amsmath}
\usepackage{graphics}
\usepackage{epic}
\def\cA{{\cal A}}          \def\cB{{\cal B}}          
\def\cD{{\cal D}}                    \def\cF{{\cal F}}

                    \def\cR{{\cal R}}
\def\cS{{\cal S}}                    \def\cU{{\cal U}}
          \def\cW{{\cal W}}

\def\fA{{\mathfrak A}}

\newcommand{\bb}{{\bar b}}
\newcommand{\bc}{{\bar c}}
\newcommand{\sltwo}{{(\widehat{sl(2)}_{c})}}
\newcommand{\elpa}{{{\cal A}_{q,p}\sltwo}}
\newcommand{\elpb}{{{\cal B}_{q,p,\lambda}\sltwo}}
\newcommand{\un}{\mbox{1\hspace{-1mm}I}}
\newcommand{\uq}[1]{{{\cal U}_{#1}\sltwo}}
\newcommand{\uqf}[1]{{{\cal U}^F_{#1}\sltwo}}
\newcommand{\uqv}[1]{{{\cal U}^V_{#1}\sltwo}}
\newcommand{\sln}{{(\widehat{sl(N)}_{c})}}
\newcommand{\dy}[2]{{{\cal D}Y_{#1}^{#2}(sl(2))_{c}}}

\newcommand{\sfrac}[2]{{\textstyle{\frac{#1}{#2}}}}
\newcommand{\half}{{\sfrac{1}{2}}}

\newcommand{\hyperg}[4]{{_{2}\phi_{1}\left(\begin{array}{c}{#1}\quad{#2}\\
                        {#3}\end{array};{#4}\right)}}

\def\sn{\mathop{\rm sn}\nolimits}
\def\snh{\mathop{\rm snh}\nolimits}
\def\cotan{\mathop{\rm cotan}\nolimits}

\def\QTHA{Quasitriangular Hopf Algebra (QTHA)\def\QTHA{QTHA}}
\def\QTQHA{Quasitriangular Quasi-Hopf Algebra (QTQHA)\def\QTQHA{QTQHA}}

\begin{document}

\pagestyle{empty}

\begin{center}

{\Large \textbf{Cladistics\footnote{
Cladistics :
a system of biological taxonomy that defines taxa uniquely by shared
characteristics not found in ancestral groups and uses inferred
evolutionary relationships to arrange taxa in a branching hierarchy such
that all members of a given taxon have the same ancestors.}
of Double Yangians and Elliptic Algebras}}


\vspace{10mm}

{\large D. Arnaudon$^a$, J. Avan$^b$, L. Frappat$^a$, E. Ragoucy$^a$, M.
Rossi$^c$}

\vspace{10mm}

\emph{$^a$ Laboratoire d'Annecy-le-Vieux de Physique Th{\'e}orique LAPTH}

\emph{CNRS, UMR 5108, associ{\'e}e {\`a} l'Universit{\'e} de Savoie}

\emph{LAPP, BP 110, F-74941 Annecy-le-Vieux Cedex, France}

\vspace{7mm}

\emph{$^b$ LPTHE, CNRS, UMR 7589, Universit{\'e}s Paris VI/VII, France}

\vspace{7mm}

\emph{$^c$ Department of Mathematics, University of Durham \\
South Road, Durham DH1 3LE, UK}

\end{center}

\vfill
\vfill

\begin{abstract}
  A self-contained description of algebraic structures, obtained by
  combinations of various limit procedures applied to vertex and face
  $sl(2)$ elliptic quantum affine algebras, is given. New double
  Yangians structures of dynamical type are in particular
  defined. Connections between these structures are established.
  A number of them take the form of twist-like actions. These
  are conjectured to be evaluations of universal twists.
\end{abstract}

\vfill
\vfill

\rightline{LAPTH-738/99}
\rightline{PAR-LPTHE 99-23}
\rightline{DTP-99-45}
\rightline{math.QA/9906189}
\rightline{June 1999}

\cleardoublepage

\baselineskip=16pt
\tableofcontents

\newpage
\pagestyle{plain}
\setcounter{page}{1}
\baselineskip=16pt


\section{Introduction}

\subsection{Overview}

The study of elliptic quantum algebras, defined with the help of
elliptic $R$-matrices, has yielded a number of algebraic structures
relevant to certain integrable systems in quantum mechanics and
statistical mechanics (noticeably the $XYZ$ model \cite{JKKMW}, RSOS models
\cite{ABF,Ko1} and Sine--Gordon theory \cite{Ko2,KLCP}).
More recently the definition and construction of some scaling limits
has led to the notion of deformed double Yangian algebras. We will
investigate and develop here in great detail the occurrence of these
and other limit algebraic structures and the pattern of connection in
between, in the simplest case of an underlying $sl(2)$ algebra.

Two classes of elliptic solutions to the Yang--Baxter equation have
been identified, respectively associated with the vertex statistical
models \cite{Ba,Be} and the face-type statistical models
\cite{ABF,DJMO,JMO}. The vertex elliptic $R$-matrix for $sl(2)$ was
first used by Sklyanin \cite{Skl} to construct a two-parameter
deformation of the enveloping algebra $\cU(sl(2))$.
The central extension of this structure was proposed in \cite{FIJKMY}
for $sl(2)$, and later extended to
${\cal A}_{q,p}\sln$ in
\cite{JKOS}. Its connection to $q$-deformed Virasoro and $\cW_N$
algebras \cite{AKOS,FR,FF} was established in \cite{AFRS3,AFRS5}.

The face-type $R$-matrices, depending on the extra parameters
$\lambda$ belonging to the dual of the Cartan algebra in the
underlying algebra, were first used by Felder \cite{Fe} to define the
algebra $\elpb$ in the $R$-matrix
approach. Enriquez
and Felder \cite{EF} and Konno \cite{Ko1} introduced a current
representation, although differences arise in the treatment of the
central extension.
A slightly different structure, also based upon face-type $R$-matrices
but incorporating extra, Heisenberg algebra generators, was introduced
as $\cU_{q,p}(sl(2))$ \cite{Ko1,JKOS2}. This structure is
relevant to the resolution of the quantum Calogero--Moser and
Ruijsenaar--Schneider models \cite{ABB,BBB,JurcoSchupp}.
Another dynamical elliptic algebra, denoted ${\cal A}_{q,p;\pi}\sltwo$,
was also defined and studied in \cite{HouYang}.  It was then interpreted, at
the level of representation, as a twist of $\elpa$.

Particular limits of the $\cA_{q,p}$-type algebras were subsequently
defined and compared with previously known structures. The limit
$p\rightarrow 0$ together with the renormalization of the generators
by suitable powers of $p$ \emph{before} taking the limit, leads to
the quantum algebra $\uq{q}$ such as presented in
\cite{RSTS,FF}. It differs from the presentation in \cite{DF} by a
scalar factor in the $R$-matrix.
The scaling limit of the algebra ${\cal A}_{q,p;\pi}\sltwo$ was also
defined in \cite{HouYang}.

A second limit was considered in \cite{JKM,Ko2} ($R$-matrix
formulation) and \cite{KLP} (current algebra formulation). It is defined by
taking $p=q^{2r}$ (elliptic nome) and $z=q^{i\beta/\pi}$ (spectral
parameter) with $q\rightarrow 1$. This algebra, denoted
$\cA_{\hbar,\eta}\sltwo$, where $\eta\equiv \frac1r$ and
$q\simeq e^{\epsilon \hbar}$ with $\epsilon\rightarrow 0$, is relevant
to the study of the $XXZ$
model in its gapless regime \cite{JKM}. It admits a further limit
$r\rightarrow \infty$ ($\eta \rightarrow 0$) where its $R$-matrix
becomes identical to the $R$-matrix defining the double Yangian
$\dy{}{}$ (centrally extended), defined in \cite{KT} (Yangian double),
\cite{Kh} (central extension); alternative versions with a different
normalization are given in \cite{IK} (for $sl(2)$) and \cite{Io} (for
$sl(N)$). This difference in the normalization factors of the
$R$-matrix, crucial in confronting the centrally extended versions, is
the exact counterpart of the difference between the presentation of
$\uq{q}$ in \cite{DF} and \cite{RSTS}.

One must however be careful in this identification in terms of
$R$-matrix structure since the generating functionals (Lax matrices)
of these algebras admit different interpretations in terms of modes
(generators of the enveloping algebra). In the context of
$\cA_{\hbar,0}\sltwo$ the expansion is done in terms of
continuous-index Fourier modes of the spectral parameter (see
\cite{Ko2,KLP});  in the
context of $\dy{}{}$ the expansion is done in terms of powers of the
spectral parameter (see \cite{KT,Kh,Io}).

It was shown recently that both vertex algebras
$\cA_{q,p}\sln$ and face-type algebras $\cB_{q,\lambda}\sln$
were in fact Drinfel'd twists \cite{Dr} of the quantum group
$\cU_q\sln$. Originating with the proposition of
\cite{BBB} on face-type algebras, the construction of the twist
operators was undertaken in both cases by Fr{\o}nsdal \cite{Fro1,Fro2}
and finally achieved at the level of formal universal twists in
\cite{JKOS,ABRR}.  In \cite{ABRR}, the universal twist is obtained by
solving a linear equation introduced in \cite{BR}, this equation
playing a fundamental r{\^o}le for complex continuation of $6j$ symbols.
Moreover
in the case of finite (super)algebras, the convergence of the
infinite products defining the
twists was also proved in \cite{ABRR}.
This led to a formal construction of universal  $R$-matrices for the
elliptic algebras $\cA_{q,p}$ and $\cB_{q,\lambda}$, of which the
BB and ABF $4 \!\times\! 4$ matrices are respectively (spin $1/2$)
evaluation representations.

\subsection{General settings}

Our strategy is to combine in as many patterns as possible the
different limit procedures introduced previously in the literature; to
apply them to cases not already considered, in particular the face
type algebras $\elpb$; and thus to achieve as large as possible a
self-contained network of algebraic structures extending from the
elliptic quantum affine algebras to the affine Lie algebra $\uq{}$.

Before summarizing our investigations, we must first of all define
precisely the concepts which we will use
throughout this paper, so that no ambiguity arises in our statements.

We shall deal with formal algebraic structures  defined by $R$-matrix
exchange relations between formal $2 \!\times\! 2$ matrix-valued
generating functionals denoted Lax operators, using the well-known $RLL$
formalism \cite{FRT}.
Explicit $R$-matrices here are interpreted as evaluation
representations of universal objects whenever they are known to exist, or
conjectural universal objects when not. We shall not give any
precise definition of the individual generators themselves, i.e. the
specific expression of the individual generators
in terms of spectral parameter dependent Lax operators.
These definitions would eventually give rise to
the fully explicit algebraic structure. For instance we shall not
distinguish here between the double Yangian $\dy{}{}$ and the
scaled algebra $\cA_{\hbar,0}\sltwo$.
Definition of, and identification between algebraic structures will
therefore be understood at the sole level of their $R$-matrix
presentation, except in explicitly specified cases where we are able
to state relations between the full (generator-described) exchange
structures, or even the Hopf or quasi-Hopf algebraic structures.
We consider that the existence of such relations is in any case an
indication that similar connections exist at the level of universal
algebras, to be explicitly formulated once the explicit algebra
generators are defined.

Similarly we shall manipulate $R$-matrices at the level of their
evaluation representation of spin $1/2$ ($4\!\times\!4$
matrices). Only when we
shall use the term ``universal'', will it mean the abstract algebraic
object known as universal $R$-matrix. The same will apply to twist
operators connecting (quasi)-Hopf algebraic structures \cite{Dr}, and
the $R$-matrices of the algebras. We recall that a twist operator $F$
lives in the square $\cA^{\otimes 2}$ of an algebraic structure; it
connects two coproducts in $\cA$ as
$\Delta_F(\cdot) = F\Delta(\cdot) F^{-1}$,
and two universal $R$-matrices as $R_F=F^{\pi} R F^{-1}$. Its evaluation
representation acts similarly on the evaluation representation of the
universal $R$-matrices:
\begin{equation}
  R^F_{12}=F_{21} R_{12} F^{-1}_{12} \;.
  \label{eq:twist_intro}
\end{equation}

As in the previous case of identifications of algebras, we conjecture
that occurrence of a relation of this form at the level of
evaluated $R$-matrices is an indication that a similar relation
exists at the level of universal algebras.
We shall therefore denote any such relation between evaluated
$R$-matrices as a ``twist-like action'' between two algebraic
structures respectively characterized by $R$ and $R^F$, even when we
do not have explicit proof that a universal twist exists between the
universal $R$-matrices, or the respective coproduct structures.

A connection of the form (\ref{eq:twist_intro})
where $F$ will not depend on any parameter (spectral ($z$ or $\beta$,
elliptic ($p$ or $r$) or
dynamical ($w$ or $s$)) will be termed ``rigid twist action''.

We must also introduce the notion of homothetical twist-like connection,
whereby we mean the existence of an invertible matrix $F(z)$ such that
two $R$-matrices are connected by
\begin{equation}
  \widetilde{R} = f(z,p,q) F_{21}(z^{-1}) R F_{12}(z)^{-1} \;,
  \label{eq:homothetic}
\end{equation}
where $f(z,p,q)$ is a $c$-number function.
\\
At this point, we do not have an interpretation of this kind of
relation between algebraic structure. We shall come back to this point
in the conclusion.

\subsection{General properties of $R$-matrices and twists}

All evaluated $R$-matrices in this paper will obey one of the
following equations, implying the associativity of the exchange
algebra.
\begin{itemize}
\item  Yang--Baxter equation:
  \begin{equation}
    R_{12}(z) \, R_{13}(zz') \, R_{23}(z') =
    R_{23}(z') \, R_{13}(zz') \, R_{12}(z) \;,
    \label{eq:YBE}
  \end{equation}
  \begin{equation}
    R_{12}(\beta) \, R_{13}(\beta+\beta') \,
    R_{23}(\beta') = R_{23}(\beta') \,
    R_{13}(\beta+\beta') \, R_{12}(\beta) \;,
    \label{eq:YBE2}
  \end{equation}
\item  Dynamical Yang--Baxter equation:
  \begin{equation}
        R_{12}(z,\lambda+h^{(3)}) \, R_{13}(zz',\lambda) \,
        R_{23}(z',\lambda+h^{(1)}) = R_{23}(z',\lambda) \,
        R_{13}(zz',\lambda+h^{(2)}) \, R_{12}(z,\lambda) \;,
    \label{eq:DYBE}
  \end{equation}
  \begin{equation}
        R_{12}(\beta,\lambda+h^{(3)}) \,
        R_{13}(\beta+\beta',\lambda) \,
        R_{23}(\beta',\lambda+h^{(1)}) =
        R_{23}(\beta',\lambda) \,
        R_{13}(\beta+\beta',\lambda+h^{(2)}) \,
        R_{12}(\beta,\lambda) \;,
    \label{eq:DYBE2}
  \end{equation}
\end{itemize}
depending upon the multiplicative or additive nature of the spectral
parameter.

Among the algebraic structures which we consider here, some are known
to have \QTHA\  structure (for instance $\uq{q}$, $\dy{}{}$)
\cite{Dri}, and others are \QTQHA\ \cite{Dr} (for instance $\elpa$,
$\elpb$).

Their universal $R$-matrices obey the universal Yang--Baxter equation
in the first case,
\begin{equation}
  \cR_{12} \cR_{13} \cR_{23} =
  \cR_{23} \cR_{13} \cR_{12}
  \label{eq:YBEuniv}
\end{equation}
and a more complicated Yang--Baxter-type equation
in the second case, involving a cocycle $\Phi\in
\fA\otimes\fA\otimes\fA$:
\begin{equation}
  \cR_{12} \Phi_{312} \cR_{13} \Phi_{132}^{-1} \cR_{23} \Phi_{123} =
  \Phi_{321} \cR_{23} \Phi_{231}^{-1} \cR_{13} \Phi_{213} \cR_{12} \;.
  \label{eq:quasiYBE}
\end{equation}

However, in all the cases which are considered here,
the $R$-matrices, once evaluated, obey the Yang--Baxter or dynamical
Yang--Baxter equation.

We now recall the following contingent properties of evaluated
$R$-matrices.

\begin{itemize}
\item  Unitarity:
  \begin{eqnarray}
    R_{12}(z) \, R_{21}(z^{-1}) &=& 1 \;,
    \label{eq:unitarity}
    \nonumber\\
    R_{12}(\beta) \, R_{21}(-\beta) &=& 1 \;,
  \end{eqnarray}
\item  Crossing-symmetry:
  \begin{eqnarray}
    \Big(R_{12}(x)^{t_2}\Big)^{-1} &=& \Big(R_{12}(q^2x)^{-1}\Big)^{t_2} \;,
    \label{eq:crossing}
    \nonumber\\
    \left( R_{12}(\beta)^{t_2}
    \right)^{-1}
    &=&
    \left( R_{12}(\beta - 2i\pi)^{-1}
    \right)^{t_2} \;,
  \end{eqnarray}
\end{itemize}
depending upon the multiplicative or additive nature of the spectral
parameter.

The unitarity relation is not satisfied in most cases: the already known
evaluated $R$-matrices for $\elpa$, $\elpb$, $\cU_{q,\lambda}\sltwo$ only
obey the crossing relation (\ref{eq:crossing})
\cite{IIJMNT,JKOS}. We shall meet with $R$-matrices obeying unitarity
relations at the end of the paper, but we have no proof that they do
correspond to evaluations of universal objects. We shall comment on
this in the conclusion.

We have indicated that Universal Twist Operators $\cal F$ transform a
coproduct $\Delta$ into another one
$\Delta^\cF(\cdot) = \cF \Delta(\cdot) \cF^{-1}$
and the $\cR$ matrix into $\cR^\cF = \cF_{21} \cR \cF_{12}^{-1}$.
If now $(\fA,\Delta,\cR)$ defines a quasi-triangular Hopf algebra and
$\cF$ satisfies the cocycle condition
\begin{equation}
  {\cF}_{12} (\Delta \otimes id)
  {\cF} = {\cF}_{23} (id \otimes \Delta)
  {\cF} \;.
  \label{eq:cocycle}
\end{equation}
$(\fA,\Delta^\cF,{\cal R}^\cF)$
defines again a quasi-triangular Hopf algebra. If however $\cF$
satisfies no particular cocycle-like relation,
$(\fA,\Delta^\cF,{\cR}^\cF)$   defines a \QTQHA:
$\cR^\cF$ satisfies then the YB-type equation (\ref{eq:quasiYBE}).
An interesting intermediate structure arises when $\cF$
satisfies a so-called shifted cocycle condition, depending upon a
parameter $\lambda$ such that \cite{BBB,Fe}:
\begin{equation}
  {\cF}_{12}(\lambda) (\Delta \otimes id) {\cF} =
  {\cF}_{23}(\lambda+h^{(1)}) (id \otimes \Delta) {\cF}
  \label{eq:cocycle_decale}
\end{equation}
where $h\in\fA$. In this case, $\cR^\cF$ satisfies the dynamical
Yang--Baxter equation (\ref{eq:DYBE}).

\subsection{Summary}

Our paper is divided into two parts.

\begin{figure}[htbp]
  \begin{center}
    \vspace{65mm}
    \includegraphics{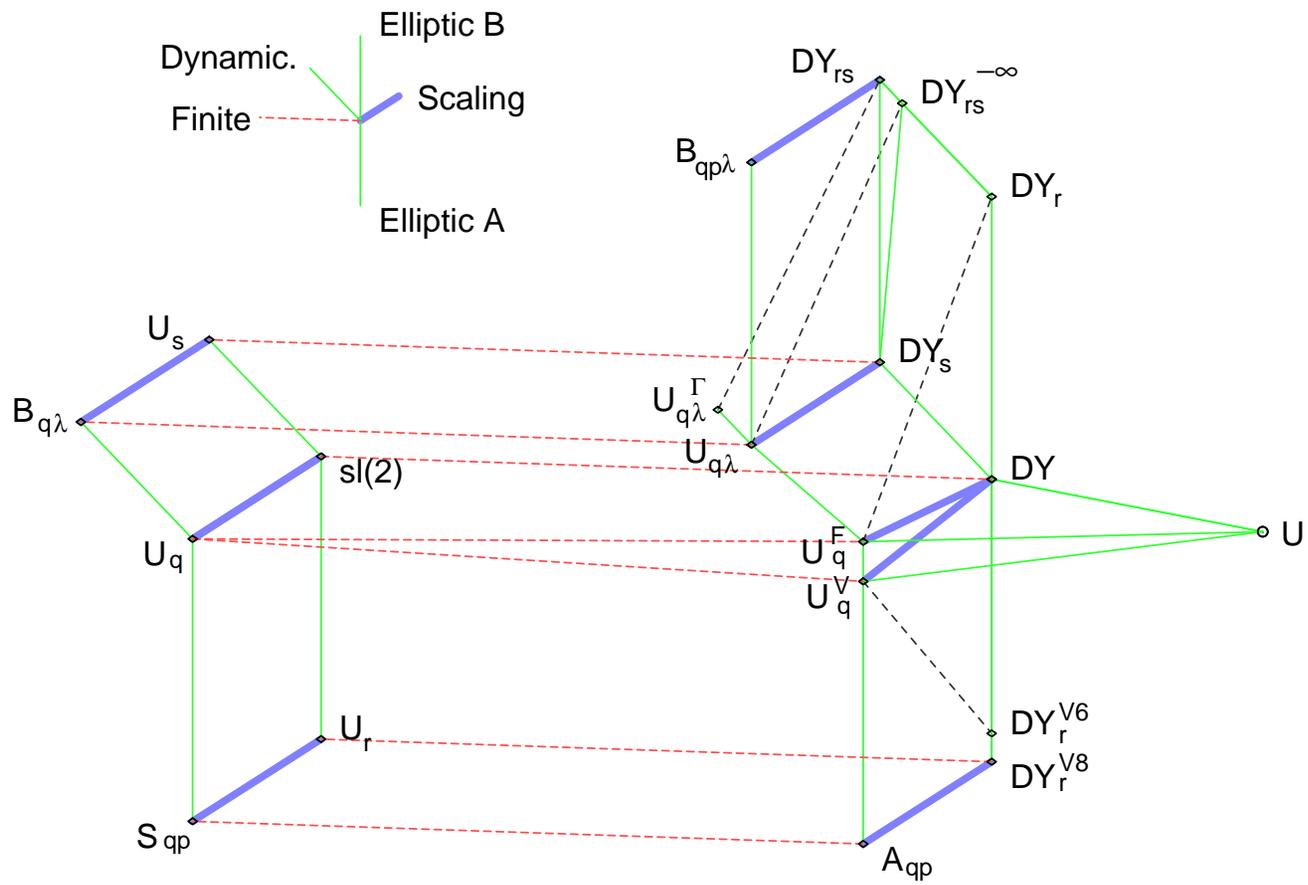}
  \end{center}
  \caption{$R$-matrix network}
  \label{fig:RCube}
\end{figure}

We shall first of all describe the limit procedures whereby the number
of parameters in the $R$-matrix description of the algebra (hence
including the spectral parameter) is decreased, starting from either
$\elpa$ or
$\elpb$; we shall define the limit
algebraic structures in both cases. These limit procedures may go in
three (for $\elpa$) or four (for $\elpb$) directions:
\begin{itemize}
\item \emph{non elliptic limit:} one sends $p$ to 0;
\item \emph{scaling limit:} one sends $q$ to 1, with $p=q^{2r}$,
  $z=q^{i\beta/\pi}$ ($z=q^{2i\beta/\pi}$ and
  $w=q^{2s}$ in the face case, where $w$ is related to $\lambda$, see
  below);
\item \emph{factorization:} one ``eliminates'' the spectral parameter by a
  Sklyanin-type
  factorization. At the level of the universal algebra this
  corresponds to a degeneracy homomorphism (see \cite{KLP}).
  This procedure is only known for vertex algebras at this
  point. Finite face type algebras however are known and shall be
  considered here, albeit without an established connection with the
  affine structures.
\item \emph{non dynamical limit:}
  in the face case the dynamical parameter $\lambda$ can also be
  eliminated by a procedure which we shall detail in the main body of
  the text.
\end{itemize}
These limit procedures, and combinations thereof, lead to
the set of objects described by Figure~\ref{fig:RCube}.
Already known structures are of course present in the diagram: $\elpb$
is the face
elliptic, centrally extended algebra; $\elpa$
is the vertex elliptic, centrally extended algebra;
$\cU_q^F\sltwo$ and $\cU_q^V\sltwo$ are two presentations
\cite{IIJMNT,JKOS} of the
quantum group $\uq{q}$ \cite{RSTS} connected by a conjugation and a
twist-like action; $\cD Y_r^{V8}\sltwo$ is the deformed double Yangian
algebra
$\cA_{\hbar,\eta}\sltwo$ in \cite{KLP} with $\hbar=1$ and
$\eta=1/r$; $\cD Y_r^{V6}\sltwo$ is the deformed double Yangian algebra
defined in \cite{Ko2}, connected to the previous one by a
rigid twist; $\dy{}{}$ is the double Yangian defined in \cite{KT,Kh};
$\cU_q(sl(2))$ is the $q$-deformed $sl(2)$ algebra; $\cS_{q,p}(sl(2))$
is Sklyanin's
elliptic ``degenerate'' algebra, and $\cU_r(sl(2))$ is the
``degenerate'' trigonometric algebra identified with $\cU_q(sl(2))$ by
$q=e^{i\pi/r}$.

New algebraic structures also appear in this diagram, mostly due to
the systematic application of the limit procedures to the face
algebra $\elpb$:  $\dy{r,s}{}$ is
the scaling limit of $\elpb$;
$\dy{r,s}{-\infty}$  is its $s\ll 0$ limit where the periodic
behaviour in $s$ is nevertheless retained;
$\dy{s}{}$ is a dynamical deformation of the double Yangian;
$\uq{q,\lambda}$ and $\cU_{q,\lambda}^\Gamma\sltwo $ are dynamical
deformations of $\uq{q}$, respectively homothetical to
$\dy{r,s}{-\infty}$ and $\dy{r,s}{}$
by a suitable redefinition of the parameters;
$\dy{r}{F}$ is an ``elliptic'' non dynamical deformation of the double
Yangian, connected to $\dy{r}{V}$ by a twist-like action and
homothetical to $\uq{q}$ by the same redefinition of the parameters.
Finally
$\cU_s(sl(2))$ and $\cB_{q,\lambda}(sl(2))$ are dynamical deformations
of the factorized structures {\`a} la Sklyanin, although they
themselves are not yet understood as originating from such a
factorization.
In addition, we also compare the structures resulting from $\elpb$ and the
structures derived \cite{HouYang} in the analysis of ${\cal
A}_{q,p;\pi}\sltwo$.  These structures are in in fact connected by a TLA
which we shall describe.

\medskip
\noindent
In order to avoid fastidious repetitions in the body of the text, we
state immediately that \emph{all} these new $R$-matrices have been
explicitly checked  to
obey the Yang--Baxter equation (\ref{eq:YBE})--(\ref{eq:YBE2}) or
dynamical
Yang--Baxter equation (\ref{eq:DYBE})--(\ref{eq:DYBE2}). Such checks
are indeed
required since the computational procedures which yield these
$R$-matrices may entail regularizations of infinite products. This
fact in
turn potentially invalidates a direct application of these
computational procedures to the Yang--Baxter equation originally
satisfied by the elliptic $R$-matrices.

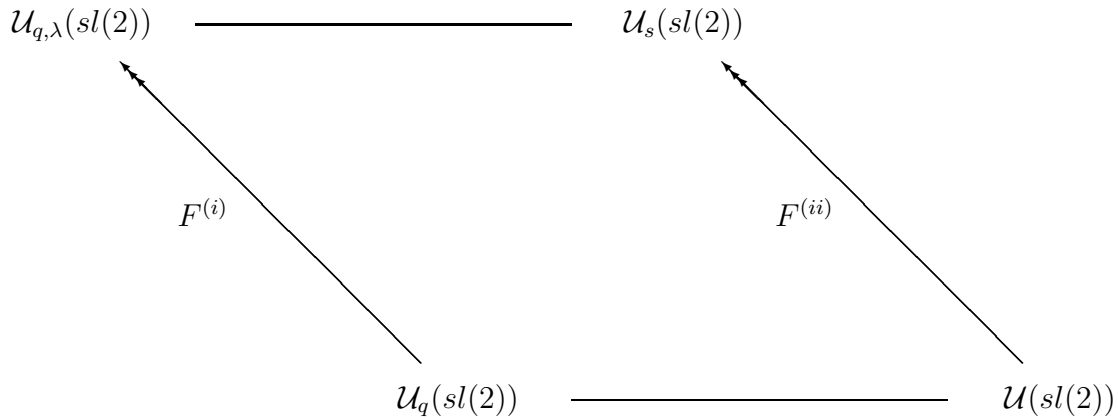
\begin{figure}[ht]
  \begin{center}
    \unitlength=1mm
    \begin{picture}(120,60)
      \put(65,50){\line(-1,0){50}}
      \put(115,0){\line(-1,0){50}}
      \put(45,5){\vector(-1,1){40}}
      \put(45,5){\vector(-1,1){39}}
      \put(45,5){\vector(-1,1){38}}
      \put(125,5){\vector(-1,1){40}}
      \put(125,5){\vector(-1,1){39}}
      \put(125,5){\vector(-1,1){38}}
      \put(50,0){\makebox(0,0){$\cU_q(sl(2))$}}
      \put(130,0){\makebox(0,0){$\cU(sl(2))$}}
      \put(0,50){\makebox(0,0){$\cU_{q,\lambda}(sl(2))$}}
      \put(80,50){\makebox(0,0){$\cU_s(sl(2))$}}
      \put(16,25){\makebox(0,0){$F^{(i)}$}}
      \put(96,25){\makebox(0,0){$F^{(ii)}$}}
    \end{picture}
  \end{center}
  \caption{The finite face case diagram: twist procedures.}
  \label{fig:ftwf}
\end{figure}

\medskip
\noindent
In the second part we describe the connections
which implement the \emph{addition} of
supplementary parameters. To be precise:
\begin{itemize}
\item implementation of the elliptic nome $p$ (or $r$);
\item implementation of the dynamical parameter $w$ (or $s$);
\item implementation of the quantum parameter $q$ along the scaling
  limit connections.
\end{itemize}
Three types of twist-like actions (TLA) appear:
\begin{enumerate}
\renewcommand{\theenumi}{\textit{\roman{enumi}}}
\item
  TLA
  explicitely proved to be evaluation of universal twists, represented
  on the figures \ref{fig:twv}--\ref{fig:ftwf} by a triple arrow. Most
  of them have  been previously
  established in the literature, particularly in
  \cite{Fro1,Fro2,Bab,JKOS}.
\item
  TLA conjectured to be evaluations of universal twists, represented
  on the figures \ref{fig:twv}--\ref{fig:ftwf}  by a double arrow. All
  these objects are new. They are
  either deduced from previously known ones by limit procedures or
  combinations; or explicitly computed from scratch.
\item
  Homothetical TLA. These are also new; they connect
  either the affine Lie algebra $\uq{}$ with double Yangian or $\uq{q}$; or
  they act as reciprocal of the scaling transformations on the vertex or
  face side.
  By contrast, let us point out that
  the first two implementations (of $p$ and $w$ -- or $r$ and $s$) are
  achieved in all cases by twist-like  actions.
\end{enumerate}
\renewcommand{\theenumi}{\arabic{enumi}}

We finally give some indications on further possible investigations in
the conclusion.

\part{Structures and Limits}

\section{Vertex type algebras}
\setcounter{equation}{0}

We will start from the elliptic algebra $\elpa$ and take the
above described different limits to obtain various quantum algebras and
deformed double Yangians.

\subsection{Elliptic algebra $\elpa$}

Let us consider the following $R$-matrix \cite{Ba,FIJKMY}:
\begin{equation}
  R(z,q,p) = \frac{\tau(q^{1/2}z^{-1})}{\mu(z)} \left(
  \begin{array}{cccc}
    a(u) & 0 & 0 & d(u) \\
    0 & b(u) & c(u) & 0 \\
    0 & c(u) & b(u) & 0 \\
    d(u) & 0 & 0 & a(u) \\
  \end{array} \right)
  \label{eq:elpa}
\end{equation}
where
\begin{eqnarray}
  && a(u) = \frac{\snh(v-u)}{\snh(v)} = z^{-1} \;
  \frac{\Theta_{p^2}(q^2z^2) \; \Theta_{p^2}(pq^2)}
  {\Theta_{p^2}(pq^2z^2) \; \Theta_{p^2}(q^2)}  \;,
  \label{eq:az}
  \\
  && b(u) = \frac{\snh(u)}{\snh(v)} = q z^{-1} \;
  \frac{\Theta_{p^2}(z^2) \; \Theta_{p^2}(pq^2)}
  {\Theta_{p^2}(pz^2) \; \Theta_{p^2}(q^2)} \;, \\
  && c(u) = 1 \;, \\
  && d(u) = -k\,\snh(v-u)\snh(u) = -p^{1/2} q^{-1} z^{-2} \;
  \frac{\Theta_{p^2}(z^2) \; \Theta_{p^2}(q^2z^2)}
  {\Theta_{p^2}(pz^2) \; \Theta_{p^2}(pq^2z^2)} \;.
  \label{eq:dz}
\end{eqnarray}
The function $\snh(u)$ is defined by $\snh(u) = -i\sn(iu)$ where $\sn(u)$ is
Jacobi's elliptic function with modulus $k$.  The variables $z,q,p$ are
related to the variables $u,v$ by
\begin{equation}
  p = \exp\Big(-\frac{\pi K'}{K}\Big) \;, \qquad
  q = - \exp\Big(-\frac{\pi v}{2K}\Big) \;, \qquad
  z = \exp\Big(\frac{\pi u}{2K}\Big) \;,
  \label{eq:pqz}
\end{equation}
where the elliptic integrals $K,K'$ are given by (with ${k'}^2=1-k^2$):
\begin{equation}
  K = \int_0^1 \frac{dx}{\sqrt{(1-x^2)(1-k^2x^2)}} \qquad\mbox{and}\qquad
  K' = \int_0^1 \frac{dx}{\sqrt{(1-x^2)(1-{k'}^2x^2)}} \;.
\end{equation}
{}From now on, we shall consider $a$, $b$, $c$, $d$, as functions of
$z$ given by (\ref{eq:pqz}).
\\
The normalization factors are
\begin{eqnarray}
  && \frac{1}{\mu(z)} = \frac{1}{\kappa(z^2)} \; \frac{(p^2;p^2)_\infty}
  {(p;p)_\infty^2} \; \frac{\Theta_{p^2}(pz^2) \; \Theta_{p^2}(q^2)}
  {\Theta_{p^2}(q^2z^2)} \;, \\
  && \frac{1}{\kappa(z^2)} = \frac{(q^4z^{-2};p,q^4)_\infty \;
    (q^2z^2;p,q^4)_\infty \; (pz^{-2};p,q^4)_\infty \; (pq^2z^2;p,q^4)_\infty}
  {(q^4z^2;p,q^4)_\infty \; (q^2z^{-2};p,q^4)_\infty \; (pz^2;p,q^4)_\infty \;
    (pq^2z^{-2};p,q^4)_\infty} \;, \\
  && \tau(q^{1/2}z^{-1}) = q^{-1/2} z \;
  \frac{\Theta_{q^4}(q^2z^2)}{\Theta_{q^4}(z^2)}  \;,
\end{eqnarray}
where the infinite multiple products are defined by:
\begin{equation}
  (z;p_1,\dots,p_m)_\infty = \prod_{n_i \ge 0} (1-z p_1^{n_1}
  \dots p_m^{n_m}) \;.
  \label{eq:prodinf}
\end{equation}
$R$ satisfies the so-called
quasi-periodicity property
\begin{equation}
  R_{12}(-z p^{\frac 12}) = (\sigma_1
  \otimes \un)^{-1} \, R_{21}(z^{-1})^{-1} \,
  (\sigma_1 \otimes \un) \;.
  \label{eq:quasiper}
\end{equation}
It also obeys the crossing-symmetry property (\ref{eq:crossing}),
but not unitarity (\ref{eq:unitarity}).
\\
This matrix defines the elliptic algebra $\elpa$ as
\begin{equation}
  R_{12}(z_1/z_2 ,q,p) \, L_1(z_1) \, L_2(z_2 ) =
  L_2(z_2 ) \, L_1(z_1) \, R_{12}(z_1/z_2 ,q,p^*=pq^{-2c}) \;.
  \label{eq:rll_elpa}
\end{equation}

\subsection{Non elliptic limit: quantum affine algebra $\uq{q}$}

Starting from the above $R$-matrix of $\elpa$, and taking the
limit $p\rightarrow 0$, one gets the $\uq{q}$ algebra, with its
$R$-matrix given by
\begin{equation}
  R_{V}(z) = \rho(z^2) \left( \begin{array}{cccc}
  1 & 0 & 0 & 0 \\
  0 & \displaystyle \frac{q(1-z^2)}{1-q^2z^2} & \displaystyle
  \frac{z(1-q^2)}{1-q^2z^2} & 0 \\
  0 & \displaystyle \frac{z(1-q^2)}{1-q^2z^2} & \displaystyle
  \frac{q(1-z^2)}{1-q^2z^2} & 0 \\
  0 & 0 & 0 & 1 \\
\end{array} \right) \;.
\label{eq:ruqa}
\end{equation}
The normalization factor is
\begin{equation}
  \rho(z^2) = q^{-1/2} \; \frac{(q^2z^2;q^4)_{\infty}^2}
  {(z^2;q^4)_{\infty} \; (q^4z^2;q^4)_{\infty}} \;.
\end{equation}
It is  known \cite{FIJKMY} that the algebra $\uq{q}$ is only
obtained after a suitable renormalization of the generators of $\elpa$
and a subsequent non-continuous limit $p \to 0$.

\noindent
The algebra $\uq{q}$ is then defined by the relations
\begin{eqnarray}
  R_{12}(z_1/z_2) \, L^\pm_1(z_1) \, L^\pm_2(z_2)
  &=&
  L^\pm_2(z_2) \, L^\pm_1(z_1) \, R_{12}(z_1/z_2) \;,
  \label{eq:rll_uq}
  \\
  R_{12}(q^{c/2} z_1/z_2) \, L^+_1(z_1) \, L^-_2(z_2)
  &=&
  L^-_2(z_2) \, L^+_1(z_1) \, R_{12}(q^{-c/2} z_1/z_2) \;.
  \label{eq:rll_uq2}
\end{eqnarray}

As indicated in the introduction, we do not discuss the problem of
generator expansions here. The same caveat
will hold throughout the whole paper, viz. we shall assume that
suitable, consistent expansions of the Lax equations will exist to
generate well-defined algebraic structures.

\subsection{Scaling limit}

The so-called scaling limit of an algebra will be understood as the
algebra defined by the scaling limit of the $R$-matrix of the initial
structure. It is obtained by setting in the $R$-matrix
$p=q^{2r}$ (elliptic nome) and $z=q^{i\beta/\pi}$ (spectral
parameter) with $q\rightarrow 1$, and $r$, $\beta$ being kept fixed.
The spectral parameter in the Lax operator is now to be taken as
$\beta$.

\subsubsection{Deformed double Yangian $\dy{r}{V8}$}
Taking the scaling limit of $\elpa$, one gets the $\dy{r}{V8}$ algebra.
Its $R$-matrix takes the form \cite{Ko2,AAFR}
(the superscript $V8$ is a token of the eight non vanishing entries of the
vertex-type $R$-matrix):
\begin{equation}
  R_{V8}(\beta,r) = \rho_{V8}(\beta;r) \left( \begin{array}{cccc}
  \displaystyle
  \frac{\cos\frac{i\beta}{2r} \; \cos\frac{\pi}{2r}}
  {\cos\frac{\pi+i\beta}{2r}}
  & 0 & 0 & \displaystyle -\frac{\sin\frac{i\beta}{2r}\;\sin\frac{\pi}{2r}}
  {\cos\frac{\pi+i\beta}{2r}} \\
  0 & \displaystyle \frac{\sin\frac{i\beta}{2r}\;\cos\frac{\pi}{2r}}
  {\sin\frac{\pi+i\beta}{2r}} & \displaystyle \frac{\cos\frac{i\beta}{2r} \;
    \sin\frac{\pi}{2r}} {\sin\frac{\pi+i\beta}{2r}} & 0 \\
  0 & \displaystyle \frac{\cos\frac{i\beta}{2r}\;\sin\frac{\pi}{2r}}
  {\sin\frac{\pi+i\beta}{2r}} & \displaystyle \frac{\sin\frac{i\beta}{2r} \;
    \cos\frac{\pi}{2r}} {\sin\frac{\pi+i\beta}{2r}} & 0 \\
  \displaystyle -\frac{\sin\frac{i\beta}{2r} \; \sin\frac{\pi}{2r}}
  {\cos\frac{\pi+i\beta}{2r}} & 0 & 0 & \displaystyle
  \frac{\cos\frac{i\beta}{2r}\;\cos\frac{\pi}{2r}}
  {\cos\frac{\pi+i\beta}{2r}} \\
\end{array} \right) \;.
\label{eq:dyra8}
\end{equation}
The normalization factor is
\begin{equation}
  \rho_{V8}(\beta;r) = -\frac{S_{2}(-\frac{i\beta}{\pi} \;\vert\; r,2) \,
    S_{2}(1+\frac{i\beta}{\pi} \;\vert\; r,2)}
  {S_{2}(\frac{i\beta}{\pi} \;\vert\; r,2) \,
    S_{2}(1-\frac{i\beta}{\pi} \;\vert\; r,2)}
  \cotan \frac{i\beta}{2} \;.
  \label{eq:rhoV8}
\end{equation}
$S_{2}(x \vert \omega_{1},\omega_{2})$ is the Barnes' double sine
function of periods $\omega_{1}$ and $\omega_{2}$ defined by
\cite{Barnes}, quoted in \cite{JimMiw}:
\begin{equation}
  S_{2}(x \vert \omega_{1},\omega_{2}) =
  \frac{\Gamma_{2}(\omega_{1}+\omega_{2}-x \;\vert\; \omega_{1},\omega_{2})}
  {\Gamma_{2}(x \;\vert\; \omega_{1},\omega_{2})}
\end{equation}
where $\Gamma_r$ is the multiple Gamma function of order $r$
given by
\begin{equation}
  \Gamma_{r}(x \vert \omega_{1},\dots,\omega_{r}) = \exp
  \left(
    \frac{\partial}{\partial s} \sum_{n_{1},\dots,n_{r} \ge 0}
    (x+n_{1}\omega_{1}+\dots+n_{r}\omega_{r})^{-s}\Bigg\vert_{s=0}
  \right) \;.
\end{equation}
This ${R}$ matrix satisfies the
quasi-periodicity property
\begin{equation}
  {R}_{12}(\beta-i\pi r) = (\sigma_1 \otimes \un)^{-1} \,
  {R}_{21}(-\beta)^{-1} \, (\sigma_1 \otimes \un) \;,
  \label{eq:quasiper2}
\end{equation}
where $\sigma_1$ is the usual Pauli matrix.
\\
It also obeys the crossing-symmetry property (\ref{eq:crossing}),
but not (\ref{eq:unitarity}).

\noindent
The algebra $\dy{r}{V8}$ is then defined by the relation
\begin{equation}
  R_{12}(\beta_{1}-\beta_{2},r) \, L_1(\beta_{1}) \, L_2(\beta_{2}) =
  L_2(\beta_{2}) \, L_1(\beta_{1}) \,
  R_{12}(\beta_{1}-\beta_{2},r-c) \;.
  \label{eq:rll_dy8}
\end{equation}

\subsubsection{Double Yangian $\dy{}{}$}

Starting now from the quantum affine algebra $\uq{q}$ and taking its
scaling limit, one obtains the double Yangian algebra
$\dy{}{}$ \cite{KT}. Its $R$-matrix is given by

\begin{equation}
  R(\beta) = \rho(\beta) \left( \begin{array}{cccc}
  1 & 0 & 0 & 0 \\
  0 & \displaystyle \frac{i\beta}{i\beta+\pi} & \displaystyle
  \frac{\pi}{i\beta+\pi} & 0 \\[.3cm]
  0 & \displaystyle \frac{\pi}{i\beta+\pi} & \displaystyle
  \frac{i\beta}{i\beta+\pi} & 0 \\
  0 & 0 & 0 & 1 \\
\end{array} \right) \;.
\label{eq:dy}
\end{equation}
The normalization factor is
\begin{equation}
  \rho(\beta) = \frac{\Gamma_{1}(\frac{i\beta}{\pi} \;\vert\; 2) \;
    \Gamma_{1}(2+\frac{i\beta}{\pi} \;\vert\; 2)}
  {\Gamma_{1}(1+\frac{i\beta}{\pi} \;\vert\; 2)^2} \;.
\end{equation}
Taking the limit $r\rightarrow \infty$ of the $R$-matrix of $\dy{r}{V8}$
(corresponding to the previous $p\rightarrow 0$ limit), one also gets
the double Yangian algebra.

\medskip
\noindent
Notice that in both previous cases, the limit procedure may be applied
directly to the Lax matrices, leading to the explicit, continuous
labelled algebras, respectively denoted ${{\cal
    A}_{\hbar,\eta}\sltwo}$ and ${{\cal A}_{\hbar,0}\sltwo}$
\cite{KLP}.

\bigskip

\begin{figure}[ht]
  \begin{center}
    \unitlength=1mm
    \begin{picture}(80,60)
      \put(15,50){\vector(1,0){50}}
      \put(15,0){\vector(1,0){50}}
      \put(0,10){\vector(0,1){30}}
      \put(80,10){\vector(0,1){30}}
      \put(0,50){\makebox(0,0){$\uq{q}$}}
      \put(0,0){\makebox(0,0){$\elpa$}}
      \put(80,50){\makebox(0,0){$\dy{}{}$}}
      \put(80,0){\makebox(0,0){$\dy{r}{V8}$}}
      \put(40,53){\makebox(0,0){scaling $q \rightarrow 1$}}
      \put(40,47){\makebox(0,0){$z=q^{i\beta/\pi}$}}
      \put(40,3){\makebox(0,0){scaling $q \rightarrow 1$}}
      \put(40,-3){\makebox(0,0){$z=q^{i\beta/\pi},p=q^{2r}$}}
      \put(-10,25){\makebox(0,0){$p \rightarrow 0$}}
      \put(90,25){\makebox(0,0){$r \rightarrow \infty$}}
    \end{picture}
  \end{center}
  \caption{The vertex case diagram: limit procedures}
  \label{fig:lmv}
\end{figure}
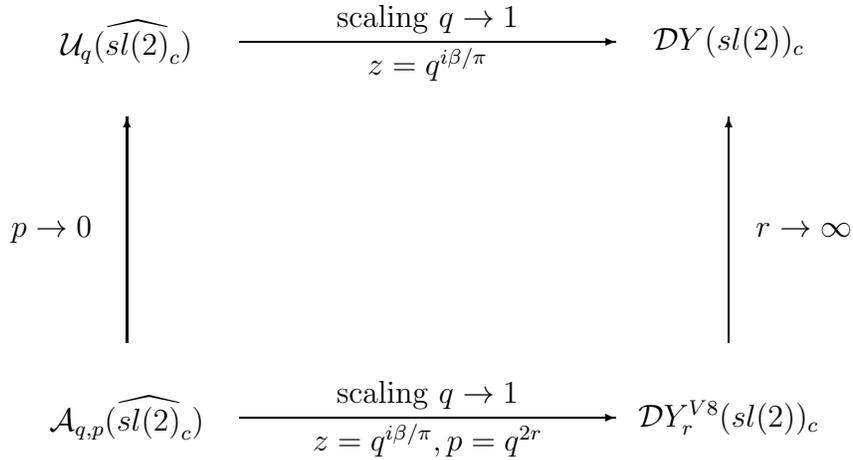

\noindent The different limit procedures in the vertex case are
summarized in Figure \ref{fig:lmv}.

\subsection{Finite algebras }
Up to now, the various limits led to affine structures. We now
consider another kind of limit where the algebra is ``factorized''. The
resulting structure is based on a finite $sl(2)$ algebra.
This is interpreted as a highly degenerate consistent representation
of the affine algebras at $c=0$, where all generators are expressed in
terms of only four ones.

\subsubsection{Sklyanin algebra}
The Sklyanin algebra \cite{Skl} is constructed from $\elpa$ taken at
$c=0$. The $R$-matrix (\ref{eq:elpa}) can be
written as
\begin{equation}
  R(z) = \un \otimes \un + \sum_{\alpha =1}^3 W_\alpha(z) \sigma_\alpha
  \otimes \sigma_\alpha \;,
  \label{eq:RSklyanin}
\end{equation}
where $\sigma_\alpha$ are the Pauli matrices and $W_\alpha(z)$ are
expressed in terms of the Jacobi elliptic functions.
A particular $z$-dependence
of the $L(z)$
operators is chosen, leading to a factorization of the $z$-dependence
in the $RLL$ relations. Indeed, setting
\begin{equation}
  L(z)= S_0 + \sum_{\alpha=1}^3 W_\alpha(z) S_\alpha \sigma_\alpha  \;,
  \label{eq:LSklyanin}
\end{equation}
one obtains an algebra with four generators
$S^\alpha$ ($\alpha=0,...,3$) and commutation relations
\begin{eqnarray}
  [S_0, S_\alpha ] &=& -i J_{\beta \gamma } (S_\beta S_\gamma + S_\gamma
  S_\beta ) \;, \nonumber \\
  {}[S_\alpha ,S_\beta ] &=& i (S_0 S_\gamma + S_\gamma  S_0) \;,
  \label{eq:Sklyanin}
\end{eqnarray}
where $J_{\alpha \beta}=\displaystyle\frac{W_\alpha^2 -
  W_\beta^2}{W_\gamma^2 - 1}$
and $\alpha$, $\beta$, $\gamma$ are cyclic permutations of 1, 2, 3.
The structure functions $J_{\alpha \beta}$ are actually
independent of $z$. Hence we get an algebra where the $z$-dependence has
been dropped out.

\subsubsection{${\cal U}_{r}(sl(2))$}

The same factorization procedure (\ref{eq:RSklyanin}-\ref{eq:LSklyanin})
applied to $\dy{r}{V8}$ leads to a ${\cal U}_{r}(sl(2))$
algebra described by
(\ref{eq:Sklyanin}) with now
$J_{12}=-J_{31}=\tan^2 \frac{\pi}{2r}$ and $J_{23}=0$.
We recognize the algebra ${\cal U}_{q'}(sl(2))$ if we set
$q'=e^{i\pi/r}$.

\medskip

\noindent\textbf{Remark:}
The scaling limit of the Sklyanin algebra (\ref{eq:Sklyanin}) also
leads to the algebra ${\cal U}_{r}(sl(2))$.

\subsubsection{Other factorizations}

Applying the factorization procedure
(\ref{eq:RSklyanin}-\ref{eq:LSklyanin})
to the quantum affine algebra $\uq{q}$, one simply gets the finite
${\cal U}_{q}(sl(2))$  algebra.
\\
Let us remark that this algebra is also the $p\rightarrow 0$
limit of the Sklyanin algebra.

\bigskip

\noindent
If we finally apply the factorization procedure to the double Yangian
$\dy{}{}$, one gets $J_{\alpha \beta}=0$. Setting the central
generator $S_0$ to 1, we recognize the classical ${\cal U}(sl(2))$
algebra.

\medskip

\noindent
Note that ${\cal U}(sl(2))$ can also be viewed  as:

\textit{i)} the $r\rightarrow \infty$ limit of ${\cal U}_{r}(sl(2))$;

\textit{ii)} the  $q\rightarrow 1$ limit (``scaling limit'') of
${\cal U}_{q}(sl(2))$.

\begin{figure}[ht]
  \begin{center}
    \unitlength=1mm
    \begin{picture}(80,60)
      \put(15,50){\vector(1,0){50}}
      \put(15,0){\vector(1,0){50}}
      \put(0,10){\vector(0,1){30}}
      \put(80,10){\vector(0,1){30}}
      \put(0,50){\makebox(0,0){$\cU_q(sl(2))$}}
      \put(0,0){\makebox(0,0){$\cS_{qp}(sl(2))$}}
      \put(80,50){\makebox(0,0){$\cU(sl(2))$}}
      \put(80,0){\makebox(0,0){$\cU_r(sl(2))$}}
      \put(40,53){\makebox(0,0){$q \rightarrow 1$}}
      \put(40,3){\makebox(0,0){scaling $q \rightarrow 1$}}
      \put(40,-3){\makebox(0,0){$p=q^{2r}$}}
      \put(-10,25){\makebox(0,0){$p \rightarrow 0$}}
      \put(90,25){\makebox(0,0){$r \rightarrow \infty$}}
    \end{picture}
  \end{center}
  \caption{The finite vertex case diagram: limit procedures}
  \label{fig:flmv}
\end{figure}
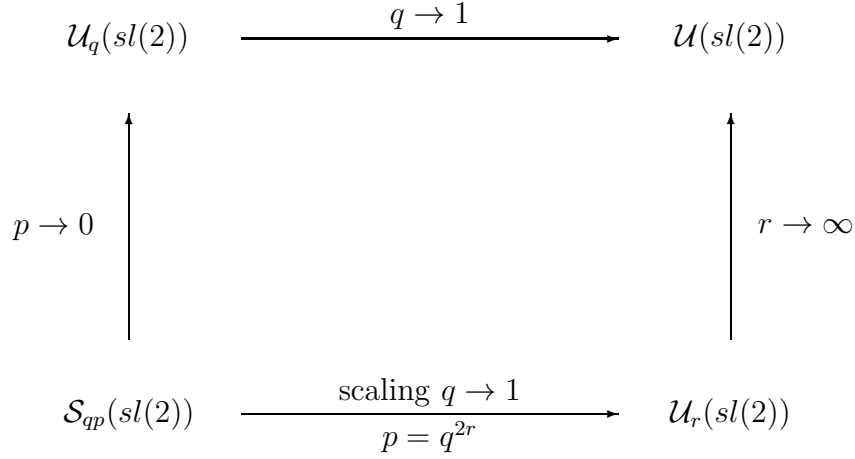

The different limit procedures in the finite vertex case are
summarized in figure \ref{fig:flmv}.

\section{Face type algebras}
\setcounter{equation}{0}

\subsection{Elliptic algebra $\elpb$}

The starting point in the face case is the $\elpb$ algebra.  Let
$\{h,c,d\}$ be a basis of the Cartan subalgebra of
$\sltwo$. If
$r,s,s'$ are complex numbers, we set $\lambda = \half\;(s+1)h + s'c +
(r+2)d$.  The elliptic parameter $p$ and the dynamical parameter $w$ are
related to the deformation parameter $q$ by $p=q^{2r}$, $w=q^{2s}$.  \\
The $R$ matrix of $\elpb$ is \cite{Fe,JKOS}
\begin{equation}
  R(z;p,w) = \rho(z;p)
  \left(
    \begin{array}{cccc}
      1 & 0 & 0 & 0 \\
      0 & b(z) & c(z) & 0 \\
      0 & \bc(z) & \bb(z) & 0 \\
      0 & 0 & 0 & 1 \\
    \end{array}
  \right)
  \label{eq:Relpb}
\end{equation}
where
\begin{eqnarray}
  b(z) &=& q \frac{(pw^{-1}q^2;p)_{\infty} \; (pw^{-1}q^{-2};p)_{\infty}}
  {(pw^{-1};p)_{\infty}^2} \; \frac{\Theta_{p}(z)}{\Theta_{p}(q^2z)}
  \;, \\
  \bb(z) &=& q \frac{(wq^2;p)_{\infty} \; (wq^{-2};p)_{\infty}}
  {(w;p)_{\infty}^2} \; \frac{\Theta_{p}(z)}{\Theta_{p}(q^2z)} \;, \\
  c(z) &=& \frac{\Theta_{p}(q^2)}{\Theta_{p}(w)} \;
  \frac{\Theta_{p}(wz)}{\Theta_{p}(q^2z)} \;, \\
  \bc(z) &=& z \frac{\Theta_{p}(q^2)}{\Theta_{p}(pw^{-1})} \;
  \frac{\Theta_{p}(pw^{-1}z)}{\Theta_{p}(q^2z)} \;.
\end{eqnarray}
The normalization factor is
\begin{equation}
  \rho(z;p) = q^{-1/2} \frac{(q^2z;p,q^4)_{\infty}^2}
  {(z;p,q^4)_{\infty} \; (q^4z;p,q^4)_{\infty}} \;
  \frac{(pz^{-1};p,q^4)_{\infty} \; (pq^4z^{-1};p,q^4)_{\infty}}
  {(pq^2z^{-1};p,q^4)_{\infty}^2} \;.
  \label{eq:rhoelpb}
\end{equation}
The elliptic algebra $\elpb$ is then defined by \cite{Fe,JKOS}
\begin{equation}
  R_{12}(z_1/z_2,\lambda+h) \, L_1(z_1,\lambda) \, L_2(z_2,\lambda+h^{(1)})
  =
  L_2(z_2,\lambda) \, L_1(z_1,\lambda+h^{(2)}) \, R_{12}(z_1/z_2,\lambda) \;.
  \label{eq:rll_elpb}
\end{equation}

\subsection{Dynamical quantum affine algebras $\uq{q,\lambda}$}

Starting from the $\elpb$ $R$-matrix, and taking the limit
$p\rightarrow 0$, one gets the $\uq{q,\lambda}$ one.
\\
The $R$ matrix of $\uq{q,\lambda}$ is
\begin{equation}
  R(z;w) = \rho(z) \left( \begin{array}{cccc}
  1 & 0 & 0 & 0 \\
  0 & \displaystyle \frac{q(1-z)}{1-q^2z} & \displaystyle
  \frac{(1-q^2)(1-wz)}{(1-q^2z)(1-w)} & 0 \\
  0 & \displaystyle \frac{(1-q^2)(z-w)}{(1-q^2z)(1-w)} & \displaystyle
  \frac{q(1-z)}{(1-q^2z)}\,\frac{(1-wq^2)(1-wq^{-2})}{(1-w)^2} & 0 \\
  0 & 0 & 0 & 1 \\
\end{array} \right) \;.
\label{eq:DQA}
\end{equation}
The normalization factor is
\begin{equation}
  \rho(z) = q^{-1/2} \; \frac{(q^2z;q^4)_{\infty}^2}{(z;q^4)_{\infty} \;
    (q^4z;q^4)_{\infty}} \;.
\end{equation}

\subsection{Non dynamical limit}

Taking the limit $w \rightarrow 0$ in $\uq{q,\lambda}$, one gets
the algebra $\uq{q}$ with $R$-matrix:

\begin{equation}
  R_{F}(z) = \rho(z) \left( \begin{array}{cccc}
  1 & 0 & 0 & 0 \\
  0 & \displaystyle \frac{q(1-z)}{1-q^2z} & \displaystyle
  \frac{1-q^2}{1-q^2z} & 0 \\
  0 & \displaystyle \frac{z(1-q^2)}{1-q^2z} & \displaystyle
  \frac{q(1-z)}{1-q^2z} & 0 \\
  0 & 0 & 0 & 1 \\
\end{array} \right) \;.
\label{eq:ruqb}
\end{equation}
The normalization factor is
\begin{equation}
  \rho(z) = q^{-1/2} \; \frac{(q^2z;q^4)_{\infty}^2}{(z;q^4)_{\infty} \;
    (q^4z;q^4)_{\infty}} \;.
  \label{eq:rhouqb}
\end{equation}

\medskip

\noindent
\textbf{Remark 1:} The matrix (\ref{eq:ruqa}) differs from the matrix
(\ref{eq:ruqb}) by rescaling $z \rightarrow z^2$ and symmetrization between
the $e_{12} \otimes e_{21}$ and $e_{21} \otimes e_{12}$ terms.
The corresponding algebraic structures will be denoted respectively
$\uqf{q}$ for (\ref{eq:ruqb}) and $\uqv{q}$ for (\ref{eq:ruqa}).

\medskip

\noindent
Actually, the matrix $R(z)$ is computed from the universal ${\cal R}$
matrix of $\uq{q}$ by
$R(z) = (\pi \otimes \pi) {\cal R}(z)$ where $\pi$ is a spin $1/2$ evaluation
representation \cite{IIJMNT}. Implementation of the spectral
parameter $z$ in the universal ${\cal R}$ matrix is obtained by
\begin{eqnarray}
  {\cal R}(z) &=& Ad(z^{\rho} \otimes 1) {\cal R} \qquad \mbox{in the
  vertex case,} \\
  {\cal R}(z) &=& Ad(z^d \otimes 1) {\cal R} \qquad \mbox{in the face case.}
\end{eqnarray}
Hence, the $R$ matrix of $\uqv{q}$ is associated to the principal gradation
of the $\sltwo$ algebra, whilst the $R$ matrix $\uqf{q}$ is
associated to the homogeneous gradation.

\medskip

\noindent
\textbf{Remark 2:} The scaling limit of the $R$-matrix (\ref{eq:ruqb})
of $\uqf{q}$ gives back the $R$-matrix (\ref{eq:dy}) of $\dy{}{}$.

\subsection{Dynamical deformed double Yangian $\dy{r,s}{}$}

Starting again from the $\elpb$ case, and taking now the scaling limit
$p=q^{2r}$ (elliptic nome), $z=q^{2i\beta/\pi}$ (spectral
parameter), $w=q^{2s}$ (dynamical parameter) with $q\rightarrow 1$,
one gets a new structure, interpreted as
a dynamical deformed centrally extended double Yangian
$\dy{r,s}{}$.
\\
The $R$ matrix of $\dy{r,s}{}$ is
\begin{equation}
  R(\beta;r,s) = \rho(\beta;r) \left( \begin{array}{cccc}
  1 & 0 & 0 & 0 \\
  0 & b(\beta) & c(\beta) & 0 \\
  0 & \bc(\beta) & \bb(\beta) & 0 \\
  0 & 0 & 0 & 1 \\
\end{array} \right) \;,
\label{eq:dyrs}
\end{equation}
where
\begin{eqnarray}
  b(\beta) &=& \frac{\Gamma_{1}(r-s \;\vert\; r)^2}
  {\Gamma_{1}(r-s+1 \;\vert\; r) \, \Gamma_{1}(r-s-1 \;\vert\; r)} \;
  \frac{\sin\frac{i\beta}{r}}{\sin\frac{\pi+i\beta}{r}} \;, \\
  c(\beta) &=& \frac{\sin\frac{\pi s+i\beta}{r}}{\sin\frac{\pi s}{r}} \;
  \frac{\sin\frac{\pi}{r}}{\sin\frac{\pi+i\beta}{r}} \;, \\
  \bb(\beta) &=& \frac{\Gamma_{1}(s \;\vert\; r)^2}
  {\Gamma_{1}(s+1 \;\vert\; r) \, \Gamma_{1}(s-1 \;\vert\; r)} \;
  \frac{\sin\frac{i\beta}{r}}{\sin\frac{\pi+i\beta}{r}} \;, \\
  \bc(\beta) &=& \frac{\sin\frac{\pi s-i\beta}{r}}{\sin\frac{\pi s}{r}} \;
  \frac{\sin\frac{\pi}{r}}{\sin\frac{\pi+i\beta}{r}} \;.
\end{eqnarray}
The normalization factor is the same as formula (\ref{eq:rhoV8}),
rewritten as
\begin{equation}
  \rho(\beta;r) = \frac{S_{2}^2(1+\frac{i\beta}{\pi} \;\vert\; r,2)}
  {S_{2}(\frac{i\beta}{\pi} \;\vert\; r,2) \,
    S_{2}(2+\frac{i\beta}{\pi} \;\vert\; r,2)} \;.
  \label{eq:normdyrs}
\end{equation}

\noindent
The algebra $\dy{r,s}{}$ is then defined by the relations
\begin{equation}
  R_{12}(\beta_1-\beta_2,\lambda+h) \,
  L_1(\beta_1,\lambda) \, L_2(\beta_2,\lambda+h^{(1)})
  =
  L_2(\beta_2,\lambda) \, L_1(\beta_1,\lambda+h^{(2)}) \,
  R_{12}(\beta_1-\beta_2,\lambda) \;.
  \label{eq:rll_dyrs}
\end{equation}

\subsection{Dynamical double Yangian $\dy{s}{}$}

Taking the limit $r\rightarrow \infty$ in $\dy{r,s}{}$, one gets
a new, dynamical, centrally extended double Yangian $\dy{s}{}$.
\\
The $R$ matrix of $\dy{s}{}$ is given by
\begin{equation}
  R(\beta) = \rho(\beta) \left( \begin{array}{cccc}
  1 & 0 & 0 & 0 \\
  0 & \displaystyle\frac{i\beta}{i\beta+\pi} &
  \displaystyle\frac{\pi s+i\beta}{s(i\beta+\pi)} & 0 \\
  0 & \displaystyle\frac{\pi s-i\beta}{s(i\beta+\pi)} &
  \displaystyle\frac{s^2-1}{s^2}\,\frac{i\beta}{i\beta+\pi} & 0 \\
  0 & 0 & 0 & 1 \\
\end{array} \right)  \;.
\label{eq:dys}
\end{equation}
The normalization factor is
\begin{equation}
  \rho(\beta) = \frac{\Gamma_{1}(\frac{i\beta}{\pi} \;\vert\; 2) \;
    \Gamma_{1}(2+\frac{i\beta}{\pi} \;\vert\; 2)}
  {\Gamma_{1}(1+\frac{i\beta}{\pi} \;\vert\; 2)^2} \;.
\end{equation}

\medskip

\noindent \textbf{Remark 1:} This $R$-matrix (\ref{eq:dys}) is also
obtained by taking the scaling limit of the $R$-matrix (\ref{eq:DQA}) of
$\uq{q,\lambda}$.

\medskip

\noindent \textbf{Remark 2:} The $|s| \rightarrow \infty$ limit of the
$R$-matrix (\ref{eq:dys}) gives back the $R$-matrix (\ref{eq:dy}) of
$\dy{}{}$.

\subsection{Dynamical deformed double Yangian $\dy{r,s}{-\infty}$ in the
  trigonometric limit \label{subsect:dyrsi}}

Starting  again from $\dy{r,s}{}$ and taking $s\ll 0$, but retaining the
oscillating behaviour in $s$, one gets $\dy{r,s}{-\infty}$, another
dynamical deformed centrally extended double Yangian structure.
\\
The $R$ matrix of $\dy{r,s}{-\infty}$ reads
\begin{equation}
  R(\beta;r,s) = \rho(\beta;r)
  \left(
    \begin{array}{cccc}
      1 & 0 & 0 & 0 \\
      0 & \displaystyle
      \frac{\sin\frac{i\beta}{r}}{\sin\frac{\pi+i\beta}{r}} &
      \displaystyle \frac{\sin\frac{\pi s+i\beta}{r}}{\sin\frac{\pi
      s}{r}} \;
      \frac{\sin\frac{\pi}{r}}{\sin\frac{\pi+i\beta}{r}} & 0 \\
      0 & \quad \displaystyle \frac{\sin\frac{\pi
      s-i\beta}{r}}{\sin\frac{\pi s}{r}}
      \; \frac{\sin\frac{\pi}{r}}{\sin\frac{\pi+i\beta}{r}} \quad &
      \quad
      \displaystyle \frac{\sin\pi\frac{s+1}{r}
      \sin\pi\frac{s-1}{r}}{\sin^2\frac{\pi
          s}{r}} \frac{\sin\frac{i\beta}{r}}{\sin\frac{\pi+i\beta}{r}}
      \quad & 0 \\
      0 & 0 & 0 & 1 \\
    \end{array}
  \right) \;.
  \label{eq:dyrsi}
\end{equation}
The normalization factor is the same as for $\dy{r,s}{}$, see
(\ref{eq:normdyrs}).

\medskip

\noindent \textbf{Remark 1:} The limit $r\rightarrow \infty$ of the
$R$-matrix (\ref{eq:dyrsi}) gives again the $R$-matrix (\ref{eq:dys})
of $\dy{s}{}$.

\medskip

\noindent \textbf{Remark 2: Correspondence with $\uq{q,\lambda}$}
\\
The previous $R$ matrix is homothetical to that of $\uq{q,\lambda}$ by
the following identifications of parameters:
\begin{equation}
  z = e^{-2\beta/r} \;, \qquad
  q = e^{i\pi/r} \;, \qquad
  w = e^{2i\pi s/r} \;.
  \label{eq:identif}
\end{equation}
The same identification of parameters applied to the $R$-matrix
(\ref{eq:dyrs}) of $\dy{r,s}{}$ leads to an $R$-matrix close to that
of $\uq{q,\lambda}$, but with $\Gamma$-function dependence in the
dynamical parameter.
This would define a new dynamical algebraic structure
${\cal U}^{\Gamma}_{q,\lambda}\sltwo$.

\subsection{Deformed double Yangian $\dy{r}{F}$}

Taking now the limit $s\rightarrow i\infty$ in $\dy{r,s}{}$, one gets
a non dynamical structure $\dy{r}{F}$.
\\
The $R$ matrix of $\dy{r}{F}$ is given by
\begin{equation}
  R(\beta;r) = \rho(\beta;r)
  \left(
    \begin{array}{cccc}
      1 & 0 & 0 & 0 \\
      0 &
      \displaystyle\frac{\sin\frac{i\beta}{r}}{\sin\frac{\pi+i\beta}{r}} &
      e^{\beta/r}
      \displaystyle\frac{\sin\frac{\pi}{r}}{\sin\frac{\pi+i\beta}{r}}
      &
      0 \\
      0 & \quad e^{-\beta/r}
      \displaystyle\frac{\sin\frac{\pi}{r}}{\sin\frac{\pi+i\beta}{r}}
      \quad & \quad
      \displaystyle\frac{\sin\frac{i\beta}{r}}{\sin\frac{\pi+i\beta}{r}}
      \quad & 0
      \\
      0 & 0 & 0 & 1 \\
    \end{array}
  \right) \;.
  \label{eq:dyrF}
\end{equation}
The normalization factor is the same as for $\dy{r,s}{}$.

\medskip

\noindent \textbf{Remark 1:} The limit $r\rightarrow \infty$ of the
$R$-matrix (\ref{eq:dyrF}) gives again the $R$-matrix of $\dy{}{}$.

\medskip

\noindent \textbf{Remark 2: Correspondence with $\uq{q}$}
\\
This matrix is homothetical to that of $\uq{q}$ -- eq.
(\ref{eq:ruqb}) -- by the following identifications of parameters:
\begin{equation}
  z = e^{-2\beta/r} \;, \qquad q = e^{i\pi/r} \;.
  \label{eq:identif2}
\end{equation}

\bigskip

\begin{figure}[ht]
  \begin{center}
    \unitlength=1mm
    \begin{picture}(120,130)
      \put(15,120){\vector(1,0){50}}
      \put(15,50){\vector(1,0){50}}
      \put(65,0){\vector(1,0){50}}
      \put(0,110){\vector(0,-1){50}}
      \put(5,45){\vector(1,-1){40}}
      \put(80,110){\vector(0,-1){50}}
      \put(85,45){\vector(1,-1){40}}
      \put(130,60){\vector(0,-1){50}}
      \put(98,114){\vector(1,-1){13}}
      \put(121,89){\vector(1,-1){11}}
      \put(128,10){\line(-1,6){13}}
      \put(10,115){\vector(1,-3){35}}
      \put(90,115){\vector(1,-3){35}}
      \put(108,90){\vector(-3,-4){22}}
      \put(50,0){\makebox(0,0){$\uq{q}$}}
      \put(130,0){\makebox(0,0){$\dy{}{}$}}
      \put(0,50){\makebox(0,0){$\uq{q,\lambda}$}}
      \put(80,50){\makebox(0,0){$\dy{s}{}$}}
      \put(0,120){\makebox(0,0){$\elpb$}}
      \put(80,120){\makebox(0,0){$\dy{r,s}{}$}}
      \put(115,95){\makebox(0,0){$\dy{rs}{-\infty}$}}
      \put(132,70){\makebox(0,0){$\dy{r}{}$}}
      \put(40,123){\makebox(0,0){scaling $q \rightarrow 1$}}
      \put(40,117){\makebox(0,0){$z=q^{2i\beta/\pi},p=q^{2r},w=q^{2s}$}}
      \put(40,53){\makebox(0,0){scaling $q \rightarrow 1$}}
      \put(40,47){\makebox(0,0){$z=q^{2i\beta/\pi},w=q^{2s}$}}
      \put(90,3){\makebox(0,0){scaling $q \rightarrow 1$}}
      \put(90,-3){\makebox(0,0){$z=q^{2i\beta/\pi}$}}
      \put(140,35){\makebox(0,0){$r \rightarrow \infty$}}
      \put(-8,85){\makebox(0,0){$p \rightarrow 0$}}
      \put(72,85){\makebox(0,0){$r \rightarrow \infty$}}
      \put(15,25){\makebox(0,0){$w \rightarrow 0$}}
      \put(95,25){\makebox(0,0){$|s| \rightarrow \infty$}}
      \put(92,67){\makebox(0,0){$r \rightarrow \infty$}}
      \put(112,110){\makebox(0,0){$s \ll 0$}}
      \put(137,85){\makebox(0,0){$s \rightarrow i\infty$}}
      \put(33,75){\makebox(0,0){$\begin{array}{c} p \rightarrow 0 \\ w
          \rightarrow 0 \end{array}$}}
      \put(117,50){\makebox(0,0){$\begin{array}{c} r \rightarrow
      \infty \\ s
          \rightarrow i\infty \end{array}$}}
      \drawline(127.3,11.5)(128,10)(128.2,11.6)(127.3,11.5)
    \end{picture}
  \end{center}
  \caption{The face case diagram: limit procedures}
  \label{fig:lmf}
\end{figure}
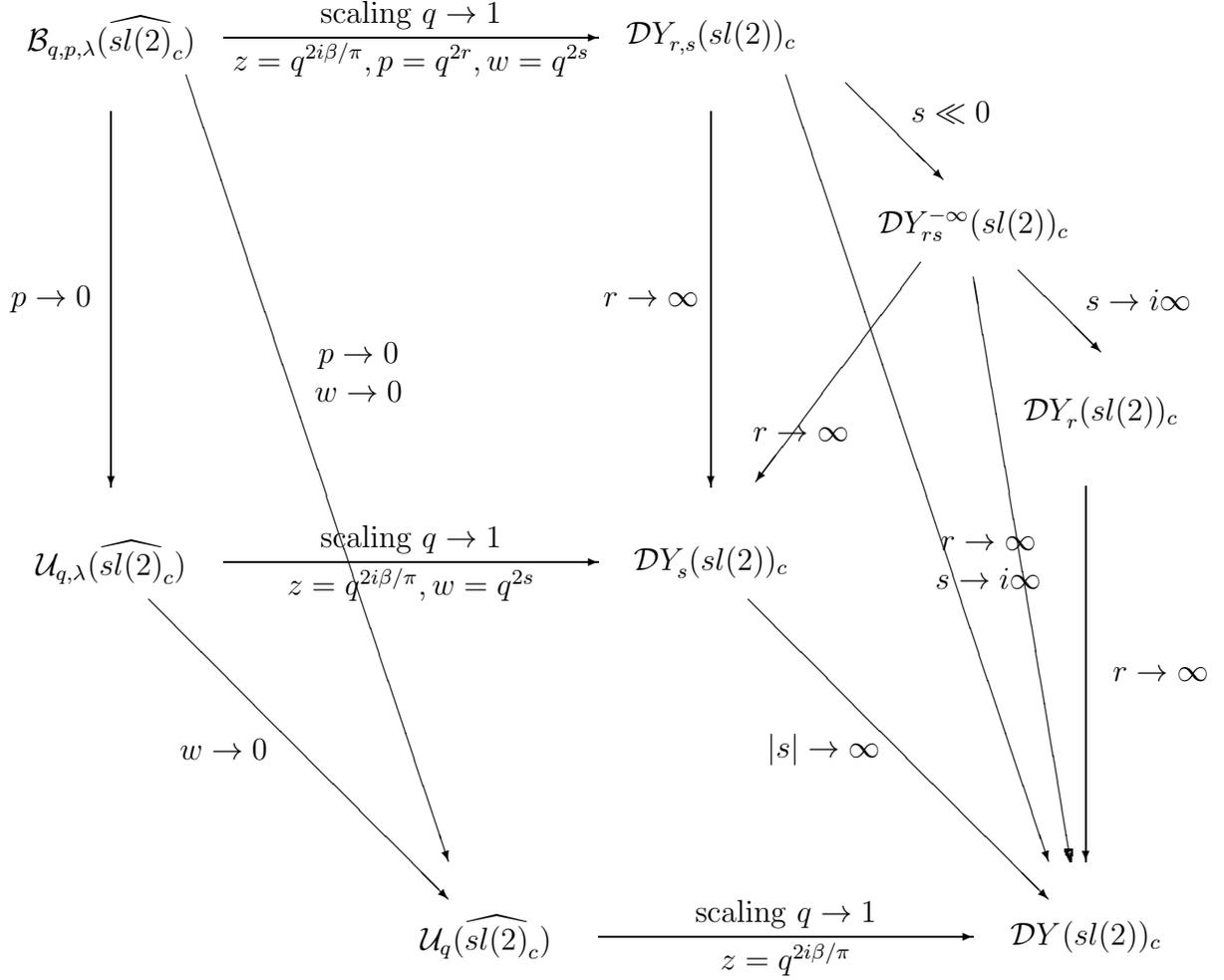

\noindent The different limit procedures in the face case are
summarized in Figure \ref{fig:lmf}.

\subsection{Finite dimensional algebras}

By constrast with the vertex case, the finite face-type elliptic
algebras have not yet been obtained from the affine algebras by a
factorization procedure. The starting point of our description will
therefore be the $R$-matrix representation of
${\cal B}_{q,\lambda}(sl(2))$ given in \cite{Bab}.

\subsubsection{Elliptic algebra ${\cal B}_{q,\lambda}(sl(2))$}

The $R$-matrix of ${\cal B}_{q,\lambda}(sl(2))$ is
\begin{equation}
  R(w) = q^{-1/2}
  \left(
    \begin{array}{cccc}
      1 & 0 & 0 & 0 \\
      0 & q & \displaystyle
      \frac{1-q^2}{1-w} & 0 \\
      0 & \displaystyle -\frac{w(1-q^2)}{1-w}\; & \displaystyle
      \;\frac{q(1-wq^2)(1-wq^{-2})}{(1-w)^2} & 0 \\
      0 & 0 & 0 & 1 \\
    \end{array}
  \right) \;.
\end{equation}

\medskip

\noindent
\textbf{Remark:} The limit $w\rightarrow 0 $ of this matrix gives the
usual $R$-matrix of $\cU_q(sl(2))$
\begin{equation}
  R =  q^{-1/2}
  \left(
    \begin{array}{cccc}
      1 & 0 & 0 & 0 \\
      0 & q & 1-q^2 & 0 \\
      0 & 0 & q & 0 \\
      0 & 0 & 0 & 1 \\
    \end{array}
  \right) \;.
  \label{eq:uqfini}
\end{equation}

\subsubsection{Dynamical algebra ${\cal U}_{s}(sl(2))$}

Taking the scaling limit $w=q^{2s}$ with $q\rightarrow 1$, we obtain
the dynamical algebra $\cU_s(sl(2))$ with the $R$-matrix
\begin{equation}
  R =
  \left(
    \begin{array}{cccc}
      1 & 0 & 0 & 0 \\
      0 & 1 & s^{-1} & 0 \\
      0 & -s^{-1} & 1-s^{-2} & 0 \\
      0 & 0 & 0 & 1 \\
    \end{array}
  \right) \;.
  \label{eq:us}
\end{equation}

The limit $|s| \rightarrow \infty$ of (\ref{eq:us}) gives $\un$, which
is the evaluated $R$-matrix of $\cU(sl(2))$.
It is not clear to us whether this particular matrix (\ref{eq:us}) can
be used for an $RLL$ formulation of the algebra. However, we will show
in the second part that it is indeed obtained as evaluation of a
\emph{universal} twist action on the \emph{universal} $R$-matrix
$1 \otimes 1$ of $\cU(sl(2))$.

\part{Twist operations}

We now describe the twist connections between the various algebraic
structures previously defined.
We first discuss twist-like actions between vertex-type algebras; we
then introduce TLAs between $\uqv{q}$ and $\uqf{q}$. We then give the
TLA between face-like algebras. The TLAs are classified here according
to the ``arrival'' algebraic structure, i.e. with the highest number
of parameters. We end up with the homothetical TLAs.

\section{Vertex type algebras}
\setcounter{equation}{0}

\subsection{Twist operator $\uqv{q} \rightarrow \elpa$}

The existence of a twist operator between $\uqv{q}$ and $\elpa$ was
proved at the level of universal matrices in \cite{JKOS}. Once the
operators are evaluated, one gets
\begin{equation}
  {R}[\elpa] = E^{(1)}_{21}(z^{-1};p) \; R[\uqv{q}] \;
  {E^{(1)}_{12}}(z;p)^{-1} \;.
\end{equation}
The twist operator $E^{(1)}(z;p)$ is given by \cite{Fro1,Fro2}
\begin{equation}
  E^{(1)}(z;p) = \rho_{E}(z;p) \left( \begin{array}{cccc}
  a_{E}(z) & 0 & 0 & d_{E}(z) \\
  0 & b_{E}(z) & c_{E}(z) & 0 \\
  0 & c_{E}(z) & b_{E}(z) & 0 \\
  d_{E}(z) & 0 & 0 & a_{E}(z) \\
\end{array} \right) \;,
\label{eq:E1}
\end{equation}
where
\begin{eqnarray}
  \label{eq:adE}
  a_{E}(z) \pm d_{E}(z) &=& \frac{(\mp p^{1/2}qz;p)_{\infty}}
  {(\mp p^{1/2}q^{-1}z;p)_{\infty}} \\
  b_{E}(z) \pm c_{E}(z) &=& \frac{(\mp pqz;p)_{\infty}}
  {(\mp pq^{-1}z;p)_{\infty}} \;.
  \label{eq:bcE}
\end{eqnarray}
The normalization factor is
\begin{equation}
  \rho_{E}(z;p) = \frac{(pz^2;p,q^4)_{\infty} \; (pq^4z^2;p,q^4)_{\infty}}
  {(pq^2z^2;p,q^4)_{\infty}^2} \;.
  \label{eq:normrhoE}
\end{equation}

\subsection{Deformed double Yangians $\dy{r}{V}$}

\subsubsection{Deformed double Yangian $\dy{r}{V6}$}

We need to define an algebraic structure not previously derived in
this paper.
\\
The $R$ matrix (\ref{eq:dyra8}) of the deformed double Yangian
$\dy{r}{V8}$  can be related to the two-body $S$ matrix of the
Sine--Gordon theory $S_{SG}(\beta,r)$ by a rigid twist operator.
The connection goes as follows. One defines the following $R$-matrix
\cite{Ko2}:
\begin{equation}
  R_{V6}(\beta,r) = \cotan(\frac{i\beta}{2}) S_{SG}(\beta,r)
  =
  \rho_{V6}(\beta;r)
  \left(
    \begin{array}{cccc}
      1 & 0 & 0 & 0 \\[1mm]
      0 & \displaystyle
      \frac{\sin\frac{i\beta}{r}}{\sin\frac{\pi+i\beta}{r}} &
      \displaystyle \frac{\sin\frac{\pi}{r}}
      {\sin\frac{\pi+i\beta}{r}} & 0 \\[5mm]
      0 & \displaystyle \frac{\sin\frac{\pi}{r}}
      {\sin\frac{\pi+i\beta}{r}} &
      \displaystyle \frac{\sin\frac{i\beta}{r}}
      {\sin\frac{\pi+i\beta}{r}} & 0 \\[5mm]
      0 & 0 & 0 & 1 \\
    \end{array}
  \right) \;,
  \label{eq:dyra6}
\end{equation}
where $\rho_{V6}(\beta;r) = \rho_{V8}(\beta;r)$, see
(\ref{eq:rhoV8}).
This $R$-matrix is assumed to define by the $RLL$ formalism an
algebraic structure denoted $\dy{r}{V6}$.
\\
One has now
\begin{equation}
  R[\dy{r}{V8}] = K_{21} \; R[\dy{r}{V6}] \; K_{12}^{-1} \;.
\end{equation}
The rigid twist operator $K$ is given by
\begin{equation}
  K = \frac{1}{2} \left( \begin{array}{rrrr}
  1 & -i & -i & -1 \\
  -1 & -i & i & -1 \\
  -1 & i & -i & -1 \\
  1 & i & i & -1 \\
\end{array} \right) \;.
\label{eq:twistcst}
\end{equation}

\medskip

\noindent \textbf{Remark 1:}
We note that
\begin{equation}
  K = V \otimes V \qquad \mbox{with} \qquad V = \frac{1}{\sqrt{2}}
  \left(
    \begin{array}{rr} 1 & -i \\ -1 & -i \\
    \end{array}
  \right) \;.
\end{equation}
This implies an isomorphism between $\dy{r}{V8}$ and $\dy{r}{V6}$
where the Lax operators are connected by $L_{V8} = V L_{V6} V^{-1}$.

\medskip

\noindent \textbf{Remark 2:}
The rigid twist leaves invariant the $R$-matrix of the undeformed
double Yangian, upon which $V$ induces an automorphism.

\medskip

\noindent \textbf{Remark 3:}
The $R$ matrix (\ref{eq:dyra6}) is homothetical
to that of $\uqv{q}$ -- eq.
(\ref{eq:ruqa}) -- by the following identifications of parameters:
\begin{equation}
  z = e^{-\beta/r} \;, \qquad q = e^{i\pi/r} \;.
  \label{eq:identiv}
\end{equation}

\medskip

\noindent \textbf{Remark 4:}
By applying the twist $K$ to the $R$-matrix (\ref{eq:ruqa})
we obtain an $R$-matrix $R[{\cal U}^{V8}_{q}\sltwo]$
with eight non-vanishing entries which may
equivalently describe $\uq{q}$ according to remark 1.
\\
Moreover, $R[{\cal U}^{V8}_{q}\sltwo]$ appears to be homothetical to
the $R$-matrix obtained
by redefining the parameters of (\ref{eq:dyra8}) according to
(\ref{eq:identiv}).

\subsubsection{Twist operator $\dy{}{} \rightarrow \dy{r}{V8}$}

The $R$-matrix of $\dy{r}{V8}$ can be obtained from the $R$-matrix of
$\dy{}{}$ by a twist-like action:
\begin{equation}
  {R}[\dy{r}{V8}] = E^{(2)}_{21}(-\beta;r) \; R[\dy{}{}] \;
  {E^{(2)}_{12}}(\beta;r)^{-1} \;.
\end{equation}
The twist operator $E^{(2)}(\beta;r)$ is the scaling limit of the
twist operator $E^{(1)}(z,p)$, see eq. (\ref{eq:E1}). It is given by
\begin{equation}
  E^{(2)}(\beta;r) = \rho_{E}(\beta;r) \left( \begin{array}{cccc}
  a_{E}(\beta) & 0 & 0 & d_{E}(\beta) \\
  0 & b_{E}(\beta) & c_{E}(\beta) & 0 \\
  0 & c_{E}(\beta) & b_{E}(\beta) & 0 \\
  d_{E}(\beta) & 0 & 0 & a_{E}(\beta) \\
\end{array} \right) \;,
\end{equation}
where
\begin{eqnarray}
  \label{eq:aEscal}
  a_{E}(\beta) + d_{E}(\beta) &=& 1 \;, \\
  a_{E}(\beta) - d_{E}(\beta) &=&
  \frac{\Gamma_{1}(r-1+\frac{i\beta}{\pi}
    \;\vert\; 2r)} {\Gamma_{1}(r+1+\frac{i\beta}{\pi} \;\vert\; 2r)}
  \;, \\
  b_{E}(\beta) + c_{E}(\beta) &=& 1 \;, \\
  b_{E}(\beta) - c_{E}(\beta) &=&
  \frac{\Gamma_{1}(2r-1+\frac{i\beta}{\pi}
    \;\vert\; 2r)} {\Gamma_{1}(2r+1+\frac{i\beta}{\pi} \;\vert\; 2r)} \;.
  \label{eq:dEscal}
\end{eqnarray}
The normalization factor is
\begin{equation}
  \rho_{E}(\beta;r) =
  \frac{\Gamma_{2}(r+1+\frac{i\beta}{\pi} \;\vert\; r,2)^2}
  {\Gamma_{2}(r+\frac{i\beta}{\pi} \;\vert\; r,2) \;
    \Gamma_{2}(r+2+\frac{i\beta}{\pi} \;\vert\; r,2)} \;.
  \label{eq:normrhoEscal}
\end{equation}

\subsubsection{Twist operator $\dy{}{} \rightarrow \dy{r}{V6}$}

Combining the previous two twist-like actions, one gets
\begin{equation}
  {R}[\dy{r}{V6}] = E^{(3)}_{21}(-\beta;r) \; R[\dy{}{}] \;
  E^{(3)}_{12}(\beta;r)^{-1} \;.
\end{equation}
The twist operator $E^{(3)}(\beta;r)$ is given by
$E^{(3)}=K^{-1}E^{(2)}$, that is
\begin{equation}
  E^{(3)}(\beta;r) = \half \; \rho_{E}(\beta;r) \left(
  \begin{array}{cccc}
    1 & -1 & -1 & 1 \\
    i(a_{E}-d_{E})(\beta) & i(b_{E}-c_{E})(\beta) & i(c_{E}-b_{E})(\beta) &
    i(d_{E}-a_{E})(\beta) \\
    i(a_{E}-d_{E})(\beta) & i(c_{E}-b_{E})(\beta) & i(b_{E}-c_{E})(\beta) &
    i(d_{E}- a_{E})(\beta) \\
    -1 & -1 & -1 & -1 \\
  \end{array} \right) \;,
  \label{eq:E6}
\end{equation}
where $a_{E},b_{E},c_{E},d_{E}$ are given by the formulae
(\ref{eq:aEscal})--(\ref{eq:dEscal})
and the normalization factor $\rho_{E}(\beta;r)$ by (\ref{eq:normrhoEscal}).

\bigskip

\noindent The different twist procedures in the vertex case are
summarized in Figure \ref{fig:twv}.

\section{Vertex to face isomorphism}
\setcounter{equation}{0}

The two $R$ matrices (\ref{eq:ruqa}) and (\ref{eq:ruqb})
can be related by a twist operator:
\begin{equation}
  R[\uqf{q}](z^2) = K^{(6)}_{21}(z^{-1}) \; R[\uqv{q}](z) \;
  {K^{(6)}_{12}(z)}^{-1} \;.
\end{equation}
The twist operator $K^{(6)}$ is given by
\begin{equation}
  K^{(6)}(z) =
  \left(
    \begin{array}{cccc}
      1 & 0 & 0 & 0 \\
      0 & z^{-1/2} & 0& 0 \\
      0 & 0 & z^{1/2} & 0 \\
      0 & 0 & 0 & 1 \\
    \end{array}
  \right)
  =
  \left(
    \begin{array}{cc}
      1&0 \\ 0&z^{1/2}
    \end{array}
  \right)
  \otimes
  \left(
  \begin{array}{cc}
    1&0 \\ 0&z^{-1/2}
  \end{array}
  \right)
  \;.
  \label{eq:K6}
\end{equation}
which acts also as a \emph{bona fide} conjugation
since $K_{21}(z^{-1})=K_{12}(z)$.
Moreover, a redefinition of the Lax operators in
(\ref{eq:rll_uq},\ref{eq:rll_uq2})
as
\begin{eqnarray}
  L_F^-(z^2) &=&
  \left(
  \begin{array}{cc}
    z^{-1/2} &0 \\ 0& z^{1/2}
  \end{array}
  \right)
  L_V^-(z)
  \left(
    \begin{array}{cc}
    z^{1/2} &0 \\ 0& z^{-1/2}
  \end{array}
  \right)
  \\
  L_F^+(z^2) &=&
  \left(
  \begin{array}{cc}
    z^{-1/2}q^{-c/4} &0 \\ 0& z^{1/2}q^{c/4}
  \end{array}
  \right)
  L_V^+(z)
  \left(
  \begin{array}{cc}
    z^{1/2}q^{-c/4} &0 \\ 0& z^{-1/2}q^{c/4}
  \end{array}
  \right)
  \label{eq:Ltransf}
\end{eqnarray}
provides a genuine algebra isomorphism between $\uqf{q}$ and $\uqv{q}$.

\section{Face type algebras}
\setcounter{equation}{0}

\subsection{Twist operator $\uqf{q} \rightarrow \uq{q,\lambda}$}

The two $R$ matrices of $\uqf{q}$  and $\uq{q,\lambda}$
can be related by a twist operator:
\begin{equation}
  R[\uq{q,\lambda}] = F^{(3)}_{21}(w) \; R[\uqf{q}] \; {F^{(3)}_{12}}(w)^{-1}
  \label{eq:twF3} \;.
\end{equation}
The twist operator $F^{(3)}(w)$ is given by
\begin{equation}
  F^{(3)}(w) =
  \left(
    \begin{array}{cccc}
      1 & 0 & 0 & 0 \\
      0 & 1 & \displaystyle\frac{w(q-q^{-1})}{1-w} & 0 \\
      0 & 0 & 1 & 0 \\
      0 & 0 & 0 & 1 \\
    \end{array}
  \right) \;.
  \label{eq:F3}
\end{equation}

\subsection{Dynamical face elliptic affine algebra $\elpb$}

\subsubsection{Twist operator $\uqf{q} \rightarrow \elpb$}

The existence of a twist operator between $\uqf{q}$ and $\elpb$ was
proved at the level of universal matrices in \cite{JKOS}. Once the
operators are evaluated, one gets
\begin{equation}
  R[\elpb] = F^{(1)}_{21}(z^{-1};p,w) \; R[\uqf{q}] \;
  {F^{(1)}_{12}}(z;p,w)^{-1} \;.
\end{equation}
The twist operator $F^{(1)}(z;p,w)$ is given by
\begin{equation}
  F^{(1)}(z;p,w) = \rho_{F}(z;p)
  \left(
    \begin{array}{cccc}
      1 & 0 & 0 & 0 \\
      0 & X_{11}(z) & X_{12}(z) & 0 \\
      0 & X_{21}(z) & X_{22}(z) & 0 \\
      0 & 0 & 0 & 1 \\
    \end{array}
  \right) \;,
\end{equation}
where
\begin{eqnarray}
  \label{eq:Xij}
  X_{11}(z) &=& \hyperg{wq^2}{q^2}{w}{p,pq^{-2}z} \;, \\
  X_{12}(z) &=& \frac{w(q-q^{-1})}{1-w} \;
  \hyperg{wq^2}{pq^2}{pw}{p,pq^{-2}z} \;, \\
  X_{21}(z) &=& z \; \frac{pw^{-1}(q-q^{-1})}{1-pw^{-1}} \;
  \hyperg{pw^{-1}q^2}{pq^2}{p^2w^{-1}}{p,pq^{-2}z} \;, \\
  X_{22}(z) &=& \hyperg{pw^{-1}q^2}{q^2}{pw^{-1}}{p,pq^{-2}z}  \;.
\end{eqnarray}
The $q$-hypergeometric function $\hyperg{q^a}{q^b}{q^c}{q,z}$ is defined by
\begin{equation}
  \hyperg {q^a}{q^b}{q^c}{q,z} = \sum_{n=0}^\infty
  \frac{(q^a;q)_{n}(q^b;q)_{n}}{(q^c;q)_{n}(q;q)_{n}} \; z^n
  \qquad \mathrm{where} \qquad
  (x;q)_{n} = \prod_{k=0}^{n-1} (1-xq^k) \;.
\end{equation}
The normalization factor is
\begin{equation}
  \rho_{F}(z;p) = \frac{(pz;p,q^4)_{\infty} \;
    (pq^4z;p,q^4)_{\infty}}
  {(pq^2z;p,q^4)_{\infty}^2}\;.
\label{eq:normfact}
\end{equation}

\subsubsection{Twist operator $\uq{q,\lambda} \rightarrow \elpb$}

Combining the last two twist-like actions, one obtains
\begin{equation}
  R[\elpb] = F^{(5)}_{21}(z^{-1};p,w) \; R[\uq{q,\lambda}] \;
  {F^{(5)}_{12}(z;p,w)}^{-1}  \;.
\end{equation}
The twist operator $F^{(5)}(z;p,w)$ is given by $F^{(5)}=F^{(1)}
{F^{(3)}}^{-1}$, that is
\begin{equation}
  F^{(5)}(z;w) =  \rho_{F}(z;p)
  \left( \begin{array}{cccc}
    1 & 0 & 0 & 0 \\
    0 & X_{11}'(z) & X_{12}'(z) & 0 \\
    0 & X_{21}'(z) & X_{22}'(z) & 0 \\
    0 & 0 & 0 & 1 \\
  \end{array} \right)  \;,
\end{equation}
where
\begin{eqnarray}
  X'_{11}(z) &=& X_{11}(z) \;, \\
  X'_{12}(z) &=& X_{12}(z) - \frac{w(q-q^{-1})}{1-w} \; X_{11}(z) \;, \\
  X'_{21}(z) &=& X_{21}(z) \;, \\
  X'_{22}(z) &=& X_{22}(z) - \frac{w(q-q^{-1})}{1-w} \; X_{21}(z) \;.
\end{eqnarray}
and $X_{ij}(z)$ are given in (\ref{eq:Xij}).
\\
The normalization factor $\rho_{F}(z;p)$ is given by (\ref{eq:normfact}).

\subsection{Dynamical double Yangian $\dy{s}{}$}

By taking the scaling limit of the connection (\ref{eq:twF3}), one
gets
\begin{equation}
  R[\dy{s}{}] = F^{(4)}_{21}(s) \; R[\dy{}{}] \;
  {F^{(4)}_{12}}(s)^{-1} \;.
\end{equation}
The twist operator $F^{(4)}(s)$ is the scaling limit of the twist
operator $F^{(3)}$, see eq. (\ref{eq:F3}). It is given by
\begin{equation}
  F^{(4)}(s) =
  \left(
    \begin{array}{cccc}
      1 & 0 & 0 & 0 \\
      0 & 1 & -s^{-1} & 0 \\
      0 & 0 & 1 & 0 \\
      0 & 0 & 0 & 1 \\
    \end{array}
  \right) \;.
\end{equation}


\subsection{Deformed double Yangian $\dy{r}{F}$}

\subsubsection{Twist operator $\dy{r}{F}$ $\rightarrow \dy{r}{V}$}

The two deformed double Yangians $\dy{r}{V}$ and $\dy{r}{F}$ obtained
from the vertex type
algebras on one hand, and from face type algebras on the other hand, are
related by twist-like actions.  One has:
\begin{equation}
  R[\dy{r}{F}](\beta) = K^{(6)}_{21}(-\beta) \; R[\dy{r}{V6}](\beta) \;
  {K^{(6)}_{12}(\beta)}^{-1} \;.
\end{equation}
The twist operator $K^{(6)}$ is actually equal to the twist operator
(\ref{eq:K6}) by setting $z=e^{-\beta/r}$:
\begin{equation}
  K^{(6)}(\beta) =
  \left(
    \begin{array}{cccc}
      1 & 0 & 0 & 0 \\
      0 & e^{\beta/2r} & 0 & 0 \\
      0 & 0 & e^{-\beta/2r} & 0 \\
      0 & 0 & 0 & 1 \\
    \end{array}
  \right) \;.
\end{equation}
Using the rigid twist operator (\ref{eq:twistcst}), one gets also:
\begin{equation}
  R[\dy{r}{F}](\beta) = K^{(8)}_{21}(-\beta) \; R[\dy{r}{V8}](\beta) \;
  {K^{(8)}_{12}(\beta)}^{-1} \;.
\end{equation}
The twist operator $K^{(8)}$ is given by $K^{(8)}=K^{(6)} K^{-1}$,
that is
\begin{equation}
  K^{(8)}(\beta) = \frac{1}{2} \left( \begin{array}{cccc}
  1 & -1 & -1 & 1 \\
  ie^{\beta/2r} & ie^{\beta/2r} & -ie^{\beta/2r} & -ie^{\beta/2r} \\
  ie^{-\beta/2r} & -ie^{-\beta/2r} & ie^{-\beta/2r} & -ie^{-\beta/2r} \\
  -1 & -1 & -1 & -1 \\
\end{array} \right) \;.
\end{equation}
This twist provides a link between the face type and
vertex type double Yangian structures.

\subsubsection{Twist operator $\dy{}{} \rightarrow \dy{r}{F}$}

The connection between $R$-matrices of $\dy{}{}$ and $\dy{r}{F}$ can
be established by three different combinations of previously
constructed twist-like actions. These three combinations of course
give by construction the
same twist operator $F^{(7)} = K^{(8)} E^{(2)} =
K^{(6)} K^{-1} E^{(2)}$. One has therefore
\begin{equation}
  R[\dy{r}{F}] = F^{(7)}_{21}(-\beta;r,) \; R[\dy{}{}] \;
  {F^{(7)}_{12}(\beta;r)}^{-1} \;.
\end{equation}
The twist operator $F^{(7)}$ is given by
\begin{equation}
  F^{(7)}(\beta;r) = \half \; \rho_{E}(\beta;r)
  \left(
    \begin{array}{cccc}
      1 & -1 & -1 & 1 \\
      i(a_{E}-d_{E})e^{\beta/2r} & i(b_{E}-c_{E})e^{\beta/2r} &
      i(c_{E}-b_{E})e^{\beta/2r} & i(d_{E}-a_{E})e^{\beta/2r} \\
      i(a_{E}-d_{E})e^{-\beta/2r} & i(c_{E}-b_{E})e^{-\beta/2r} &
      i(b_{E}-c_{E})e^{-\beta/2r} & i(d_{E}-a_{E})e^{-\beta/2r} \\
      -1 & -1 & -1 & -1 \\
    \end{array}
  \right) \;,
  \label{eq:F7}
\end{equation}
where $a_{E},b_{E},c_{E},d_{E}$ are given by
(\ref{eq:aEscal})--(\ref{eq:dEscal})
and the normalization factor $\rho_{E}(\beta;r)$ by (\ref{eq:normrhoEscal}).

\subsection{Trigonometric Dynamical deformed double Yangian
$\dy{r,s}{-\infty}$}

\subsubsection{Twist operator $\dy{r}{F} \rightarrow \dy{r,s}{-\infty}$}

The connection between $\dy{r}{F}$ and $\dy{r,s}{-\infty}$ is achieved
by the twist operator $F^{(3)}$:
\begin{equation}
  R[\dy{r,s}{-\infty}] = F^{(3)}_{21}(s) \; R[\dy{r}{F}] \;
  {F^{(3)}_{12}(s)}^{-1}  \;.
\end{equation}
The twist operator $F^{(3)}$ is actually equal to the twist operator
(\ref{eq:F3}) by setting $q=e^{i\pi/r}$ and $w=e^{2i\pi s/r}$:
\begin{equation}
  F^{(3)}(s) =
  \left(
    \begin{array}{cccc}
      1 & 0 & 0 & 0 \\
      0 & 1 & -e^{i\pi s/r} \displaystyle
      \frac{\sin\frac{\pi}{r}}{\sin\frac{\pi s}{r}} & 0 \\
      0 & 0 & 1 & 0 \\
      0 & 0 & 0 & 1 \\
    \end{array}
  \right) \;.
\end{equation}

\subsubsection{Twist operator $\dy{}{} \rightarrow \dy{r,s}{-\infty}$}

Combination of two twist-like operations yields the connection between
$\dy{}{}$ and $\dy{r,s}{-\infty}$:
\begin{equation}
  R[\dy{r,s}{-\infty}] = F^{(10)}_{21}(-\beta;r,s) \; R[\dy{}{}] \;
  {F^{(10)}_{12}(\beta;r,s)}^{-1} \;.
\end{equation}
The twist operator $F^{(10)}$ is given by $F^{(10)} = F^{(3)}
F^{(7)}$, that is
\begin{eqnarray}
  &&
  F^{(10)}(\beta;r,s) = \half \; \rho_{E}(\beta;r)
  \nonumber\\
  && \qquad
  \left( \begin{array}{cccc}
  1 & -1 & -1 & 1 \\
  i(a_{E}-d_{E})e_{+} & i(b_{E}-c_{E})e_{-} & i(c_{E}-b_{E})e_{-} &
  i(d_{E}-a_{E})e_{+} \\
  i(a_{E}-d_{E})e^{-\beta/2r} & i(c_{E}-b_{E})e^{-\beta/2r} &
  i(b_{E}-c_{E})e^{-\beta/2r} & i(d_{E}-a_{E})e^{-\beta/2r} \\
  -1 & -1 & -1 & -1 \\
\end{array} \right) \nonumber\\
\end{eqnarray}
where $e_{\pm} = e^{\beta/2r} \mp e^{i\pi s/r} e^{-\beta/2r} \frac{\sin(\pi
/r)}{\sin(\pi s/r)}$, the functions $a_{E},b_{E},c_{E},d_{E}$ are given by
the formulae (\ref{eq:aEscal})--(\ref{eq:dEscal})
and the normalization factor
$\rho_{E}(\beta;r)$ by (\ref{eq:normrhoEscal}).

\subsubsection{Twist operator $\dy{s}{} \rightarrow \dy{r,s}{-\infty}$}

Again, combination of two twist-like operations yields the connection
between $\dy{s}{}$ and $\dy{r,s}{-\infty}$:
\begin{equation}
  R[\dy{r,s}{-\infty}] = F^{(8)}_{21}(-\beta;r,s) \; R[\dy{s}{}] \;
  {F^{(8)}_{12}(\beta;r,s)}^{-1} \;.
\end{equation}
The twist operator $F^{(8)}$ is given by $F^{(8)} = F^{(10)}
{F^{(4)}}^{-1}$, that is
\begin{eqnarray}
  && F^{(8)}(\beta;r,s) = \half \; \rho_{E}(\beta;r) \nonumber \\
  && \qquad \left( \begin{array}{cccc}
  1 & -1 & -1-s^{-1} & 1 \\
  i(a_{E}-d_{E})e_{+} & i(b_{E}-c_{E})e_{-} &
  i(c_{E}-b_{E})e_{-}(1-s^{-1}) &
  i(d_{E}-a_{E})e_{+} \\
  i(a_{E}-d_{E})e^{-\beta/2r} & i(c_{E}-b_{E})e^{-\beta/2r} &
  i(b_{E}-c_{E})e^{-\beta/2r}(1-s^{-1}) &
  i(d_{E}-a_{E})e^{-\beta/2r} \\
  -1 & -1 & -1-s^{-1} & -1 \\
\end{array} \right) \nonumber \\
\label{eq:F8}
\end{eqnarray}
where $e_{\pm}$, $a_{E},b_{E},c_{E},d_{E}$ and $\rho_{E}(\beta;r)$ have the
same meaning as above.

\subsection{Dynamical deformed double Yangian $\dy{r,s}{}$}

\subsubsection{Twist operator $\dy{r,s}{-\infty} \rightarrow \dy{r,s}{}$}

The $R$-matrices of  $\dy{r,s}{-\infty}$ and $\dy{r,s}{}$ are
connected by a diagonal TLA (not depending on the spectral parameter):
\begin{equation}
  R[\dy{r,s}{}] = G_{21}(r,s) \; R[\dy{r,s}{-\infty}] \;
  {G_{12}(r,s)}^{-1} \;.
\end{equation}
The twist operator $G$ is given by
\begin{equation}
  G(r,s) = \left( \begin{array}{cccc}
  1 & 0 & 0 & 0 \\
  0 & g^{-1} & 0 & 0 \\
  0 & 0 & g & 0 \\
  0 & 0 & 0 & 1 \\
\end{array} \right)
=
  \left(
    \begin{array}{cc}
      1&0 \\ 0&g
    \end{array}
  \right)
  \otimes
  \left(
  \begin{array}{cc}
    1&0 \\ 0&g^{-1}
  \end{array}
  \right)
\;,
\label{eq:G}
\end{equation}
where
\begin{equation}
  g(r,s) = \displaystyle\frac{\Gamma_{1}(r-s \; \vert \; r)}
  {\Gamma_{1}(r-s+1 \; \vert \; r)^{1/2} \;
    \Gamma_{1}(r-s-1 \; \vert \; r)^{1/2}}\;.
  \label{eq:defg}
\end{equation}

\medskip

\noindent \textbf{Remark:}
Equivalently, $G$ expressed in terms of the parameters $q=e^{i\pi/r}$
and $w=e^{2i\pi s/r}$ realizes a TLA between $\uq{q,\lambda}$ and
${\cal U}^{\Gamma}_{q,\lambda}\sltwo$ defined in remark 2, Section
\ref{subsect:dyrsi}.

\subsubsection{Twist operator $\dy{s}{} \rightarrow \dy{r,s}{}$}

Combining the last two twists, one gets:
\begin{equation}
  R[\dy{r,s}{}] = F^{(6)}_{21}(-\beta;r,s) \; R[\dy{s}{}] \;
  {F^{(6)}_{12}(\beta;r,s)}^{-1}  \;.
\end{equation}
The twist operator $F^{(6)}$ is given by $F^{(6)} = G F^{(8)}$, that is
\begin{eqnarray}
  && F^{(6)}(\beta;r,s) = \half \; \rho_{E}(\beta;r) \nonumber \\
  &&
  \left(
    \begin{array}{cccc}
      1 & -1 & -1-s^{-1} & 1 \\
      i(a_{E}-d_{E})e_{+}g^{-1} & i(b_{E}-c_{E})e_{-}g^{-1} &
      i(c_{E}-b_{E})e_{-}(1-s^{-1})g^{-1} & i(d_{E}-a_{E})e_{+}g^{-1} \\
      i(a_{E}-d_{E})e^{-\beta/2r}g & i(c_{E}-b_{E})e^{-\beta/2r}g &
      i(b_{E}-c_{E})e^{-\beta/2r}(1-s^{-1})g &
      i(d_{E}-a_{E})e^{-\beta/2r}g \\
      -1 & -1 & -1-s^{-1} & -1 \\
    \end{array}
  \right) \nonumber \\
\end{eqnarray}
where $e_{\pm}$, $a_{E},b_{E},c_{E},d_{E}$ and $\rho_{E}(\beta;r)$ have the
same meaning as above and $g$ is given by (\ref{eq:defg}).

\subsubsection{Twist operator $\dy{}{} \rightarrow \dy{r,s}{}$}

Similarly, by a combination of previous twists, one gets:
\begin{equation}
  R[\dy{r,s}{}] = F^{(2)}_{21}(-\beta;r,s) \; R[\dy{}{}] \;
  {F^{(2)}_{12}(\beta;r,s)}^{-1} \;.
\end{equation}
The twist operator $F^{(2)}$ is given by $F^{(2)} = G F^{(10)}$, that is
\begin{eqnarray}
  && F^{(2)}(\beta;r,s) = \half \; \rho_{E}(\beta;r) \nonumber \\
  && \qquad \left( \begin{array}{cccc}
  1 & -1 & -1 & 1 \\
  i(a_{E}-d_{E})e_{+}g^{-1} & i(b_{E}-c_{E})e_{-}g^{-1}
  & i(c_{E}-b_{E})e_{-}g^{-1} &
  i(d_{E}-a_{E})e_{+}g^{-1} \\
  i(a_{E}-d_{E})e^{-\beta/2r}g & i(c_{E}-b_{E})e^{-\beta/2r}g &
  i(b_{E}-c_{E})e^{-\beta/2r}g & i(d_{E}-a_{E})e^{-\beta/2r}g \\
  -1 & -1 & -1 & -1 \\
\end{array} \right) \nonumber \\
\end{eqnarray}
where $e_{\pm}$, $a_{E},b_{E},c_{E},d_{E}$ and $\rho_{E}(\beta;r)$ have the
same meaning as above and $g$ is given by (\ref{eq:defg}).

\subsubsection{Twist operator $\dy{r}{F} \rightarrow \dy{r,s}{}$}

Finally, connection between  $\dy{r}{F}$ and $\dy{r,s}{}$ is provided
by
\begin{equation}
  R[\dy{r,s}{}] = F^{(11)}_{21}(r,s) \; R[\dy{r}{F}] \;
  {F^{(11)}_{12}(r,s)}^{-1}  \;.
\end{equation}
The twist operator $F^{(11)}$ is given by $F^{(11)} = G F^{(3)}$, that is
\begin{equation}
  F^{(11)}(r,s) =
  \left(
    \begin{array}{cccc}
      1 & 0 & 0 & 0 \\
      0 & g^{-1} & -e^{i\pi s/r} \displaystyle
      \frac{\sin\frac{\pi}{r}}{\sin\frac{\pi s}{r}} \,  g^{-1} & 0 \\
      0 & 0 & g & 0 \\
      0 & 0 & 0 & 1 \\
    \end{array}
  \right) \;,
\end{equation}
where $g$ is given by (\ref{eq:defg}).

\bigskip

\noindent The different twist procedures in the face case are
summarized in Figure \ref{fig:twf}.

\subsection{Connections with ${\cal A}_{q,p;\pi}\sltwo $ and derived
algebras}

\subsubsection{Twist $\elpb \to {\cA}_{q,p;\pi}\sltwo $}

The $R$-matrix of ${\cal A}_{q,p;\pi}\sltwo$ given in \cite{HouYang}
(actually their $R^+$-matrix) is
\begin{equation}
        R = z^{1/2r}\rho(z) \left(
        \begin{array}{cccc}
                1 & 0 & 0 & 0 \\[1mm]
                0 &
                \displaystyle\frac{\Theta_{p}(z)\Theta_{p}(q^{-2}w)}
                {\Theta_{p}(q^2z)\Theta_{p}(w)}
                &
                \displaystyle\frac{\Theta_{p}(zw)\Theta_{p}(q^2)}
                {\Theta_{p}(q^2z)\Theta_{p}(w)}
                & 0 \\[5mm]
                0 &
                \displaystyle\frac{\Theta_{p}(z^{-1}w)\Theta_{p}(q^2)}
                {\Theta_{p}(q^2z)\Theta_{p}(w)}
                &
                \displaystyle\frac{\Theta_{p}(z)\Theta_{p}(q^2 w)}
                {\Theta_{p}(q^2z)\Theta_{p}(w)}
                & 0 \\[4mm]
                0 & 0 & 0 & 1 \\
        \end{array}
        \right)
        \label{eq:raqppi}
\end{equation}
where $\rho(z)$ is the same as (\ref{eq:rhoelpb}).
This $R$-matrix is obtained from (\ref{eq:Relpb}) by exchanging
factors in $b$ and
$\bar b$ so as to reconstruct a full $\Theta$-function dependence and
correcting the $z^{1/2r}$ factor. All this can be achieved by a
factorized diagonal twist of the form of $G$ (\ref{eq:G}).

\subsubsection{Twist $\dy{r,s}{-\infty} \to {\cA}_{\hbar,\eta;\pi}\sltwo $}

Again, the $R$-matrix of the scaling limit
${\cA}_{\hbar,\eta;\pi}\sltwo$ \cite{HouYang} of
${\cA}_{q,p;\pi}\sltwo$
can be obtained from that of $\dy{r,s}{-\infty}$ (\ref{eq:dyrsi})
by a
factorized diagonal twist. It also has the form of $G$ (\ref{eq:G}),
with now $g^2 = \frac{\sin\pi(s-1)/r}{\sin\pi s/r}$.

\subsection{Finite dimensional algebras}

In both cases where TLA actions are known for non affine algebras,
they are evaluations of universal twists.

\subsubsection{Elliptic algebra ${\cal B}_{q,\lambda}(sl(2))$}

The twist operator that links ${\cal U}_{q}(sl(2))$ to
${\cal B}_{q,\lambda}(sl(2))$ is \cite{Bab}:
\begin{equation}
  F^{(i)}(w) =  \left( \begin{array}{cccc}
  1 & 0 & 0 & 0 \\
  0 & 1 & \displaystyle\frac{(q-q^{-1})w}{1-w} & 0 \\
  0 & 0 & 1 & 0 \\
  0 & 0 & 0 & 1 \\
\end{array} \right)\;.
\end{equation}
The universal form of the twist is \cite{Bab}
\begin{equation}
  \cF(w) = \sum_{n=0}^{\infty} \frac{(q^2w)^n(q-q^{-1})^n}
  {(n)_{q^{-2}}! (q^{-2} w(t^2\otimes 1); q^{-2})_n}
  (et)^n \otimes (tf)^n \;,
  \label{eq:Funivw}
\end{equation}
where
\begin{equation}
  (x;q)_n = \prod_{i = 0}^{n-1} (1-x q^{i}) \;.
  \label{eq:prodfini}
\end{equation}

\subsubsection{Dynamical algebra ${\cal U}_{s}(sl(2))$}

Its $R$ matrix can be obtained by action of the twist
\begin{equation}
  F^{(ii)}(s) =  \left( \begin{array}{cccc}
  1 & 0 & 0 & 0 \\
  0 & 1 & -s^{-1} & 0 \\
  0 & 0 & 1 & 0 \\
  0 & 0 & 0 & 1 \\
\end{array} \right)
\end{equation}
on the $R$ matrix of ${\cal U}(sl(2))$: $R=\un_{4\times 4}$.

\medskip
\noindent
We find the universal form of the twist to be
\begin{equation}
  {\cF}= \sum_{n=0}^\infty \frac{1}{n!}\left( \prod_{k=0}^{n-1}
  [(1+k-s) 1 - h]\otimes 1\right)^{-1}\ e^n\otimes f^n \;.
\end{equation}
This twist is the scaling limit of the universal twist
(\ref{eq:Funivw}).
We checked that it satisfies the shifted cocycle condition
(\ref{eq:cocycle_decale}).

\section{Homothetical twists \label{sect:homothetical}}
\setcounter{equation}{0}

We recall that homothetical TLAs connect two $R$-matrices \emph{up to
  a scalar factor}:
\begin{equation}
  \widetilde{R} = f(z,p,q) F_{21}(z^{-1}) R F_{12}(z)^{-1} \;.
\end{equation}
From now on, we shall denote such a relation by:
\begin{equation}
  \widetilde{R} \;\sim\; F_{21}(z^{-1}) R F_{12}(z)^{-1} \;.
\end{equation}
We now describe two sets of homothetical
TLAs. The first one starts from the unit evaluated $R$-matrix of
$\uq{}$ and leads to unitary $R$-matrices. The second one goes
backwards along direction of the scaling limits.

\medskip
\noindent
It is important to notice at this point that the Lie algebraic
structure of $\uq{}$ is \emph{not} described by the $RLL$ formalism
using its unit $R$-matrix (this was also the case for $\cU(sl(2))$).
In fact, the Lie algebraic structure is described by the
semi-classical $r$-matrix, i.e. the next-to-leading order of the
evaluated universal $R$-matrix of $\uq{q}$.

\subsection{Unitary matrices}

Four homothetical TLAs can be defined between the unit matrix and
the vertex quantum affine algebras. By construction, a TLA on the unit
matrix will lead to unitary $R$-matrices, while vertex quantum affine
algebras are defined by crossing-symmetry but \emph{non}-unitary
$R$-matrices.

\subsubsection{Twist operator $\uq{} \rightarrow \uqv{q}$}

\begin{equation}
  R[\uq{q}] \;\sim\; H^{(1)}_{21}(z^{-1}) \; \mbox{1\hspace{-1mm}I}
  \;  {H^{(1)}_{12}}(z)^{-1} \;.
\end{equation}
The twist operator $H^{(1)}(z)$ is given by
\begin{equation}
  H^{(1)}(z) =  \left( \begin{array}{cccc}
  1 & 0 & 0 & 0 \\
  0 & \displaystyle
  \frac{q^{1/2}z^{(-1-\epsilon)/2} -
  q^{-1/2}z^{(1+\epsilon)/2}}{qz^{-1} - q^{-1}z}
  &
  \displaystyle
  \frac{q^{1/2}z^{(-1+\epsilon)/2} -
  q^{-1/2}z^{(1-\epsilon)/2}}{qz^{-1} - q^{-1}z}
  & 0 \\[.5cm]
  0 & \displaystyle
  \frac{q^{1/2}z^{(-1+\epsilon)/2} -
  q^{-1/2}z^{(1-\epsilon)/2}}{qz^{-1} - q^{-1}z}
  &
  \displaystyle
  \frac{q^{1/2}z^{(-1-\epsilon)/2} -
  q^{-1/2}z^{(1+\epsilon)/2}}{qz^{-1} - q^{-1}z}
  & 0 \\
  0 & 0 & 0 & 1 \\
\end{array} \right) \;,
\end{equation}
where $\epsilon$ is an arbitrary non-vanishing parameter.

\subsubsection{Twist operator $\uq{} \rightarrow \elpa$}

\begin{equation}
  R[\elpa] \;\sim\; H^{(2)}_{21}(z^{-1}) \; \mbox{1\hspace{-1mm}I}  \;
  {H^{(2)}_{12}}(z)^{-1} \;.
\end{equation}
The twist operator $H^{(2)}(z)$ is given by
\begin{equation}
  H^{(2)}(z) =
  \left(
    \begin{array}{cccc}
      A & 0 & 0 & D \\
      0 & B & C & 0 \\
      0 & C & B & 0 \\
      D & 0 & 0 & A \\
    \end{array}
  \right) \;,
\end{equation}
such that
\begin{eqnarray}
  A(1/z) \pm D(1/z) &=& [ a(z)\pm d(z)][A(z)\pm D(z)] \;,\nonumber\\
  B(1/z) \pm C(1/z) &=& [ b(z)\pm c(z)][B(z)\pm C(z)] \;,
\label{eq:ABCD}
\end{eqnarray}
the functions $a$, $b$, $c$, $d$ being the entries of the $R$-matrix
of $\elpa$ (\ref{eq:elpa}).
Solutions of (\ref{eq:ABCD}), viewed as a system of functional
equations for $A$, $B$, $C$, $D$, do exist since the functions $a \pm
d$, $b \pm c$ (\ref{eq:az})-(\ref{eq:dz}) all have precisely the form
$f(z)/f(z^{-1})$.
One can choose for instance
\begin{eqnarray}
  A(z) \pm D(z) &=& (\mp q^{-1}z^{-1}p^{1/2};p)_\infty
  (\mp qzp^{1/2};p)_\infty
  \;, \\
  B(z) \pm C(z) &=& f(z)\; (\mp pq{-1}z^{-1};p)_\infty (\mp
  pqz;p)_\infty \;,
\end{eqnarray}
where
\begin{equation}
  f(z) = \frac{q^{1/2}z^{-1}-q^{-1/2}z \pm q^{1/2} \mp q^{-1/2}}
  {qz^{-1}-q^{-1}z} \;.
\end{equation}

\subsubsection{Twist operator $\uq{} \rightarrow \dy{}{}$}

\begin{equation}
  R[\dy{}{}] \;\sim\; H^{(3)}_{21}(-\beta) \; \mbox{1\hspace{-1mm}I}
  \;  {H^{(3)}_{12}}(\beta)^{-1} \;.
\end{equation}
The twist operator $H^{(3)}(\beta)$ is given by
\begin{equation}
  H^{(3)}(\beta) =
  \left(
    \begin{array}{cccc}
      1 & 0 & 0 & 0 \\
      0 & \displaystyle
      \frac{\pi - i\beta(1+\epsilon)}
      {2(\pi - i\beta)}
      &
      \displaystyle \frac{\pi  - i\beta(1-\epsilon)}
      {2(\pi - i\beta)}
      & 0 \\[.5cm]
      0 & \displaystyle
      \frac{\pi  - i\beta(1-\epsilon)} {2(\pi - i\beta)}
      &
      \displaystyle
      \frac{\pi - i\beta(1+\epsilon)}
      {2(\pi - i\beta)}
      & 0 \\
      0 & 0 & 0 & 1 \\
    \end{array}
  \right) \;.
\end{equation}

\subsubsection{Twist operator $\uq{} \rightarrow \dy{r}{V6}$}

\begin{equation}
  R[\dy{r}{V6}] \;\sim\; H^{(4)}_{21}(-\beta) \;
  \mbox{1\hspace{-1mm}I} \; {H^{(4)}_{12}}(\beta)^{-1} \;.
\end{equation}
The twist operator $H^{(4)}(\beta)$ is given by
\begin{equation}
  H^{(4)}(\beta) =
  \left(
    \begin{array}{cccc}
      1 & 0 & 0 & 0 \\
      0 & \displaystyle
      \frac{\sin \frac{\pi - i\beta - \epsilon i\beta }{2r}}
      {\sin \frac{\pi - i\beta}{r}}
      &
      \displaystyle
      \frac{\sin \frac{\pi - i\beta + \epsilon i\beta }{2r}}
      {\sin \frac{\pi - i\beta}{r}}
      & 0 \\[.5cm]
      0 & \displaystyle
      \frac{\sin \frac{\pi - i\beta + \epsilon i\beta }{2r}}
      {\sin \frac{\pi - i\beta}{r}}
      &
      \displaystyle
      \frac{\sin \frac{\pi - i\beta - \epsilon i\beta }{2r}}
      {\sin \frac{\pi - i\beta}{r}}
      & 0 \\
      0 & 0 & 0 & 1 \\
    \end{array}
  \right) \;.
\end{equation}

\subsection{Inverse scaling procedures}

\subsubsection{Inverse scaling procedure $\dy{}{}$ to $\cU_q^V\sltwo$}

Using the correspondence (\ref{eq:identiv}), the formula
(\ref{eq:E6}) gives a homothetical twist between $\dy{}{}$ and
$\cU_q^V\sltwo$, that is, the inverse of the scaling procedure
$\cU_q^V\sltwo \rightarrow \dy{}{}$:
\begin{equation}
  R[\uqv{q}](z=e^{-\beta/r},q=e^{i\pi/r}) \;\sim\;
  E^{(3)}_{21}(-\beta;r) \; R[\dy{}{}](\beta;r) \;
  E^{(3)}_{12}(\beta;r)^{-1}  \;.
\end{equation}

\subsubsection{Inverse scaling procedure $\dy{}{V8}$ to $\elpa$}

The identification between the $R$ matrices of $\uqv{q}$ and
$\dy{r}{V6}$ through the formulae (\ref{eq:identiv}) allows us to get
a homothetical twist operator between $\dy{r}{V8}$ and $\elpa$, that
is, the inverse of the scaling procedure $\elpa \rightarrow
\dy{r}{V8}$.  More precisely, one has:
\begin{equation}
  R[\elpa](z=e^{-\beta/r},q=e^{i\pi/r}) \;\sim\;
  E^{(4)}_{21}(z^{-1};p) \; R[\dy{}{V8}](\beta;r) \;
  E^{(4)}_{12}(z;p)^{-1}  \;.
\end{equation}
The twist operator $E^{(4)}(z,p)$ is given by $E^{(4)}=E^{(1)}K^{-1}$, that is
\begin{equation}
  E^{(4)}(z,p) = \half \; \rho_{E}(z;p)
  \left(
    \begin{array}{rrrr}
      (a_{E}-d_{E}) & -(a_{E}+d_{E}) & -(a_{E}+d_{E}) & (a_{E}-d_{E}) \\
      i(b_{E}+c_{E}) & i(b_{E}-c_{E}) & i(c_{E}-b_{E}) & -i(b_{E}+c_{E}) \\
      i(b_{E}+c_{E}) & i(c_{E}-b_{E}) & i(b_{E}-c_{E}) & -i(b_{E}+c_{E}) \\
      (d_{E}-a_{E}) & -(a_{E}+d_{E}) & -(a_{E}+d_{E}) & (d_{E}-a_{E}) \\
    \end{array}
  \right) \;,
  \label{eq:E8}
\end{equation}
where $a_{E},b_{E},c_{E},d_{E}$ are given by the formulae
(\ref{eq:adE},\ref{eq:bcE})
and the normalization factor $\rho_{E}(z;p)$ by (\ref{eq:normrhoE}).

\subsubsection{Inverse scaling procedure  $\dy{}{}$ to $\cU_q^F\sltwo$}

Using the correspondence (\ref{eq:identif2}), the formula
(\ref{eq:F7}) gives a homothetical twist between $\dy{}{}$ and
$\cU_q^F\sltwo$, that is, the inverse of the scaling procedure
$\cU_q^F\sltwo \rightarrow \dy{}{}$:
\begin{equation}
  R[\uqf{q}](z=e^{-2\beta/r},q=e^{i\pi/r}) \;\sim\;
  F^{(7)}_{21}(-\beta;r) \; R[\dy{}{}](\beta;r) \;
  F^{(7)}_{12}(\beta;r)^{-1}  \;.
\end{equation}

\subsubsection{Inverse scaling procedure $\dy{s}{}$ to $\uq{q,\lambda}$}

Using the correspondence (\ref{eq:identif}), the formula
(\ref{eq:F8}) gives a homothetical twist between $\dy{s}{}$ and
$\uq{q,\lambda}$, that is, the inverse of the scaling procedure
$\uq{q,\lambda} \rightarrow \dy{s}{}$:
\begin{equation}
  R[\uq{q,\lambda}](z=e^{-2\beta/r},q=e^{i\pi/r},w=e^{2i\pi s/r})
  \;\sim\; F^{(8)}_{21}(-\beta;r,s) \;
  R[\dy{s}{}](\beta;r) \; F^{(8)}_{12}(\beta;r,s)^{-1} \;.
\end{equation}

\subsubsection{Inverse scaling procedure $\dy{r,s}{}$ to $\elpb$}

The identification between the $R$ matrices of $\uq{q,\lambda}$ and
$\dy{r,s}{-\infty}$ through the formulae (\ref{eq:identif}) allows us
to get a homothetical twist operator between $\dy{r,s}{}$ and $\elpb$,
that is, the inverse of the scaling procedure $\elpb \rightarrow
\dy{r,s}{}$. One has:
\begin{equation}
  R[\elpb] \;\sim\; F^{(9)}_{21}(z;p,w) \; R[\dy{r,s}{}] \;
  {F^{(9)}_{12}(z;p,w)}^{-1} \;.
\end{equation}
The twist operator $F^{(9)}$ is given by $F^{(9)} = F^{(5)} G^{-1} =
F^{(1)} {F^{(3)}}^{-1} G^{-1}$, that is
\begin{equation}
F^{(9)}(z;p,w) =
\left(
\begin{array}{cccc}
1 & 0 & 0 & 0 \\
0 & g X_{11} & g^{-1} (X_{12} - \frac{w(q-q^{-1})}{1-w} X_{11}) & 0 \\
0 & g X_{21} & g^{-1} (X_{22} - \frac{w(q-q^{-1})}{1-w} X_{21}) & 0 \\
0 & 0 & 0 & 1 \\
\end{array}
\right) \;,
\end{equation}
where the $X_{ij}(z)$ are given by the formulae (\ref{eq:Xij}), $g$ is
given by (\ref{eq:defg}).


\begin{figure}[ht]
  \begin{center}
    \unitlength=1mm
    \begin{picture}(120,60)
      \put(65,50){\vector(-1,0){50}}
      \put(65,0){\vector(-1,0){50}}
      \put(0,40){\vector(0,-1){30}}
      \put(0,40){\vector(0,-1){29}}
      \put(0,40){\vector(0,-1){28}}
      \put(80,40){\vector(0,-1){30}}
      \put(80,40){\vector(0,-1){29}}
      \put(0,50){\makebox(0,0){$\uq{q}$}}
      \put(0,0){\makebox(0,0){$\elpa$}}
      \put(80,50){\makebox(0,0){$\dy{}{}$}}
      \put(80,0){\makebox(0,0){$\dy{r}{V8}$}}
      \put(95,50){\vector(1,-1){20}}
      \put(95,50){\vector(1,-1){19}}
      \put(115,20){\vector(-1,-1){20}}
      \put(115,20){\vector(-1,-1){19}}
      \put(120,25){\makebox(0,0){$\dy{r}{V6}$}}
      \put(40,54){\makebox(0,0){$E^{(3)}(z;p)$}}
      \put(40,-4){\makebox(0,0){$E^{(4)}(z;p)$}}
      \put(-10,25){\makebox(0,0){$E^{(1)}(z;p)$}}
      \put(70,25){\makebox(0,0){$E^{(2)}(\beta;r)$}}
      \put(110,10){\makebox(0,0){$K$}}
      \put(115,43){\makebox(0,0){$E^{(3)}(\beta;r)$}}
    \end{picture}
  \end{center}
  \caption{The vertex case diagram: twist procedures}
  \label{fig:twv}
\end{figure}
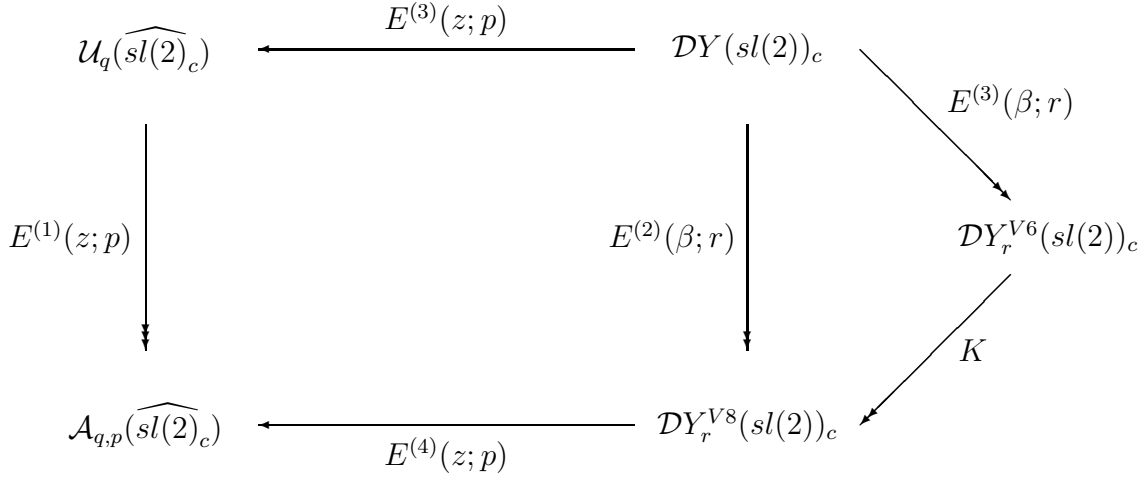

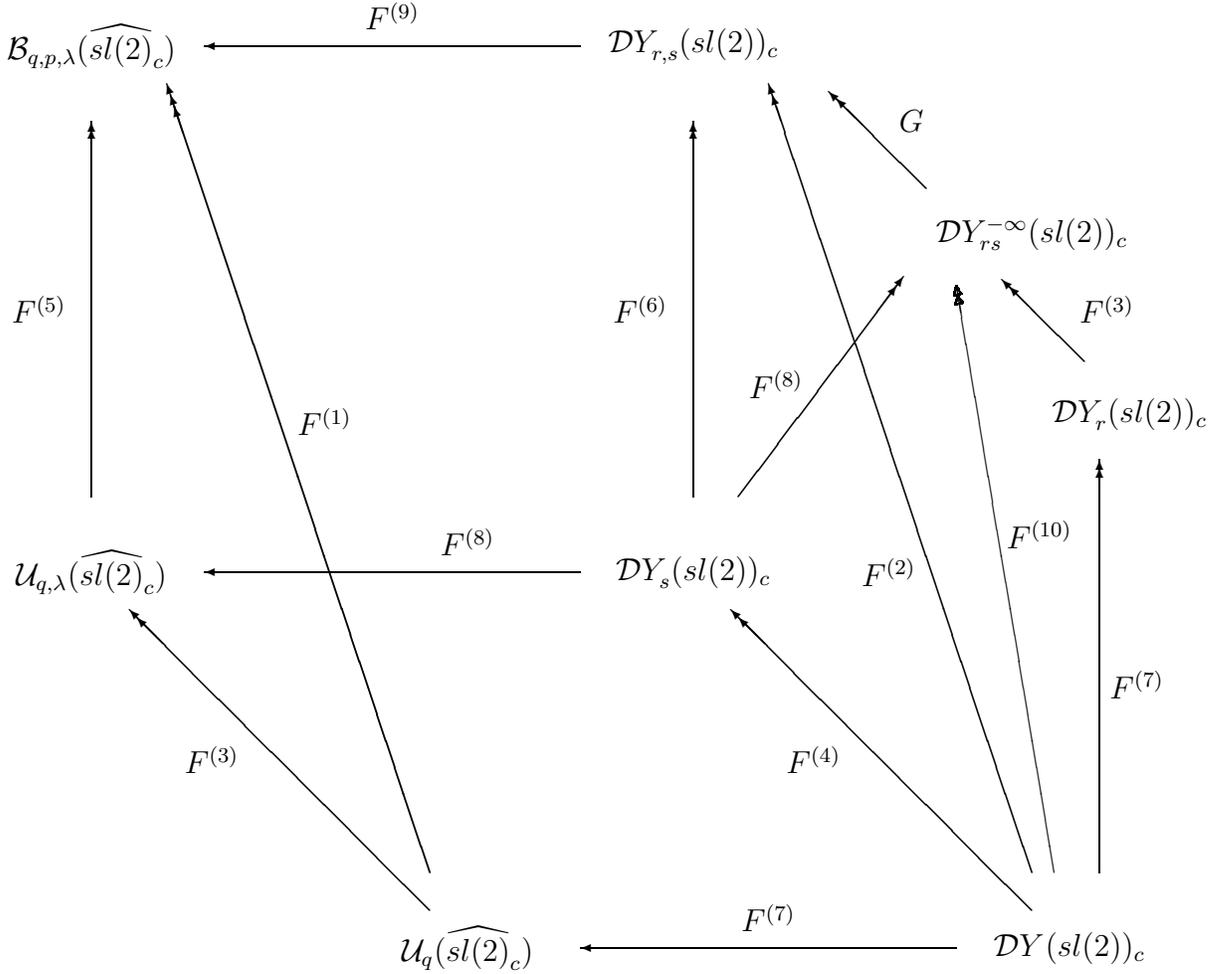
\begin{figure}[ht]
  \begin{center}
    \unitlength=1mm
    \begin{picture}(120,130)
      \put(65,120){\vector(-1,0){50}}
      \put(65,50){\vector(-1,0){50}}
      \put(115,0){\vector(-1,0){50}}
      \put(0,60){\vector(0,1){50}}
      \put(0,60){\vector(0,1){49}}
      \put(45,5){\vector(-1,1){40}}
      \put(45,5){\vector(-1,1){39}}
      \put(80,60){\vector(0,1){50}}
      \put(80,60){\vector(0,1){49}}
      \put(125,5){\vector(-1,1){40}}
      \put(125,5){\vector(-1,1){39}}
      \put(134,10){\vector(0,1){55}}
      \put(134,10){\vector(0,1){54}}
      \put(111,101){\vector(-1,1){13}}
      \put(111,101){\vector(-1,1){12}}
      \put(132,78){\vector(-1,1){11}}
      \put(132,78){\vector(-1,1){10}}
      \put(128,10){\line(-1,6){13}}
      \drawline(114.7,86.8)(115,88)(115.7,87)(114.7,86.8)
      \drawline(114.9,85.6)(115.2,86.8)(115.9,85.8)(114.9,85.6)
      \put(45,10){\vector(-1,3){35}}
      \put(45,10){\vector(-1,3){34.5}}
      \put(45,10){\vector(-1,3){34}}
      \put(125,10){\vector(-1,3){35}}
      \put(125,10){\vector(-1,3){34.5}}
      \put(86,60){\vector(3,4){22}}
      \put(86,60){\vector(3,4){21.1}}
      \put(50,0){\makebox(0,0){$\uq{q}$}}
      \put(130,0){\makebox(0,0){$\dy{}{}$}}
      \put(0,50){\makebox(0,0){$\uq{q,\lambda}$}}
      \put(80,50){\makebox(0,0){$\dy{s}{}$}}
      \put(0,120){\makebox(0,0){$\elpb$}}
      \put(80,120){\makebox(0,0){$\dy{r,s}{}$}}
      \put(125,95){\makebox(0,0){$\dy{rs}{-\infty}$}}
      \put(138,71){\makebox(0,0){$\dy{r}{}$}}
      \put(31,70){\makebox(0,0){$F^{(1)}$}}
      \put(106,50){\makebox(0,0){$F^{(2)}$}}
      \put(16,25){\makebox(0,0){$F^{(3)}$}}
      \put(96,25){\makebox(0,0){$F^{(4)}$}}
      \put(-7,85){\makebox(0,0){$F^{(5)}$}}
      \put(73,85){\makebox(0,0){$F^{(6)}$}}
      \put(90,4){\makebox(0,0){$F^{(7)}$}}
      \put(50,54){\makebox(0,0){$F^{(8)}$}}
      \put(40,124){\makebox(0,0){$F^{(9)}$}}
      \put(126,55){\makebox(0,0){$F^{(10)}$}}
      \put(135,85){\makebox(0,0){$F^{(3)}$}}
      \put(139,35){\makebox(0,0){$F^{(7)}$}}
      \put(91,75){\makebox(0,0){$F^{(8)}$}}
      \put(109,110){\makebox(0,0){$G$}}
    \end{picture}
  \end{center}
  \caption{The face case diagram: twist procedures.}
  \label{fig:twf}
\end{figure}

\clearpage

\section{Conclusion}
\setcounter{equation}{0}

We have now constructed several $R$-matrix representations for
algebraic structures, deduced from vertex or face elliptic quantum
$sl(2)$ algebras by suitable limit procedures. We have shown that
these structures exhibited associativity properties characterized by
(dynamical) Yang--Baxter equations for their evaluated
$R$-matrices. Finally, we have constructed a reciprocal set of
twist-like transformations, acting on the evaluated $R$-matrices
canonically as $R^{T}_{12} = T_{21} R_{12} T_{12}^{-1}$.

The next step is now to try to get explicit universal formulae for
these $R$-matrices and twist operators. This in turn requires to
specify the exact form under which individual generators are
encapsulated in the Lax matrices, and obtain thus the full description
of the associative algebras which we wish to study.

Let us immediately indicate that we need in particular to separate (as
is explained in \cite{KLP}) the two algebraic structures contained in
the single $R$-matrix formulations labelled here as (deformed)
(dynamical) double Yangians $\dy{...}{}$. Expansion of the Lax matrix
in terms of integer labelled generators will lead to the (deformed)
(dynamical) versions of the genuine double Yangian
\cite{KT,Kh,IK,Io}; expansion in terms of Fourier modes by a contour
integral will lead to the ``scaled elliptic'' algebras \cite{Ko2,KLP}
more correctly labelled $\cA_{\hbar,\eta}\sltwo$. Once this is done,
we can then start to investigate the following issues
\begin{itemize}
\item[--] Representations, vertex operators.
\item[--] Hopf or quasi-Hopf algebra structure, leading to:
\item[--] Universal $R$-matrices and twists.
\end{itemize}

Concerning these last two points a number of already known explicit
results lead us to draw reasonable conjectures on some of the newly
discovered algebraic structures in our work.

\subsection{Known universal $R$-matrices and twists}
Universal $R$-matrices are known for $\uq{q}$ \cite{KT2}; $\elpa$
\cite{JKOS} and $\elpb$ \cite{JKOS,Fe};
$\cB_{q,\lambda}(sl(2))$ \cite{Bab}. They are also known
for the double Yangian $\cD Y(sl(N))_c$ \cite{Kh} (proved for $N=2$,
conjectured for $N\ge 3$). Universal twists have been constructed in
the finite-algebra case from $\cU_q(sl(N))$ to
$\cB_{q,\lambda}(sl(N))$ \cite{BBB,JKOS,ABRR}; and in the affine case
from $\uq{q}$ to $\elpa$ \cite{JKOS} and $\elpb$ \cite{Fro1,Fro2}.

\subsection{Conjectures}
We therefore expect that universal $R$-matrices and twist operators
may be obtained for the complete set of algebraic structures
represented by Figure~\ref{fig:twv} in the vertex case and
Figure~\ref{fig:twf} in the face case. The structures $\dy{...}{}$ are
here to be interpreted as genuine, integer-labelled double
Yangians. The explicit construction of universal objects in this frame
seems achievable, along the lines followed in \cite{Kh} and
\cite{KT}.
The problem of constructing universal objects associated with the
continuous-labelled algebras of $\cA_{\hbar,\eta}$-type is more
delicate, since one needs in particular to contrive a direct
universal connection between continuous-labelled generators in
$\cA_{\hbar,\eta}$ and discrete-labelled generators in $\cA_{qp}$, or
between $\cA_{\hbar,0}$ and $\uq{q}$.

\subsection{The case of unitary matrices}
We have described in Section~\ref{sect:homothetical} homothetical
twist-like connections
between $\un_{4\times 4}$, interpreted as the evaluated $R$-matrix
$\un$ for the centrally extended algebra $\uq{}$, and unitary
$R$-matrices realizing a $RLL$-structure ``proportional'' to
$\uq{q}$.
Interpretation of this $RLL$-structure, and its derived relations at
elliptic level, remains obscure. The canonical construction of
universal $R$-matrices for $\uq{q}$ \cite{KT2} and their subsequent
evaluation \cite{IIJMNT} leaves open the possibility of an alternative
construction of universal $\rightarrow$ evaluated $R$-matrix which lead
to unitary (and crossing-symmetrical) $R$-matrices; it may
arise either by dropping the triangularity requirement $\cR\in
\cB_+\otimes\cB_-\subset \uq{q}\otimes\uq{q}$, or by
relaxing analyticity constraints on the evaluated $R$-matrix.
\\
Homothetical TLAs also appear between double Yangian-like structures and
their antecedent structures through the scaling procedure.  The same
possibilities hold for the differently normalized $R$-matrix structures
obtained by applications of these homothetical TLAs.

\subsection{The notion of dynamical elliptic algebra}

Finally let us briefly comment on the notion of ``dynamical'' algebraic
structure.  This notion was applied throughout this paper to algebras
incorporating an extra parameter $\lambda$ belonging to the Cartan
algebra, subsequently shifted along a general
Cartan algebra direction.  This shift is therefore retained in the
Yang-Baxter equation for evaluated $R$-matrices of face type (but
\emph{not} of vertex type, for which the extra parameter is simply a
$c$-number and the shift takes place along the central
charge direction, set to zero in the evaluation representation
\footnote{This fact was clarified to us by O. Babelon.}).  A particular
illustration of this fact arises in the case of classical and quantum
$R$-matrix for Calogero--Moser models \cite{ABB} where $\lambda$ is
identified with the momentum of the Calogero--Moser particles, hence the
denomination ``dynamical'' for the $R$-matrices.  In the algebraic
structures described here
however, $\lambda$ is not yet promoted to the r{\^o}le of generator, hence this
denomination is slightly abusive.  There exists however at least one
example of algebraic structure, $\cU_{q,p}\sltwo$ \cite{Ko1,JKOS2}, where
$\lambda$ and its conjugate $\displaystyle\frac{\partial}{\partial\lambda}$
are ``added'' to the algebra $\elpb$; however $\cU_{q,p}\sltwo$ is not
a Hopf, even quasi-Hopf, algebra.  We expect therefore that similar
genuinely dynamical algebraic structures may be associated in the same
way to all ``dynamical'' algebras described here, and may play important
r{\^o}le in solving the models where such algebras arise.

\bigskip

\noindent \textbf{Acknowledgements}

This work was supported in part by CNRS and EC network contract number
FMRX-CT96-0012.  M.R. was supported by an EPSRC research grant
no. GR/K 79437.
We wish to thank O.~Babelon for enlightening discussions,
S.M.~Khoroshkin for fruitful explanations and P.~Sorba for his helpful
and clarifying comments.
J.A.  wishes to thank the LAPTH for its kind hospitality.

\end{document}